\documentclass[final,3p,times]{elsarticle}
\usepackage{xcolor}  
\usepackage[latin1]{inputenc}
\usepackage{comment}  
\usepackage{amsfonts,amssymb,amsmath,exscale} 
\usepackage{textcomp}
\usepackage{pifont}
\usepackage{rotating}
\usepackage{wasysym}
\usepackage{defcml_fonts}
\usepackage{hyperref}
\usepackage{subcaption}

\usepackage{graphicx,dblfloatfix}
\usepackage{color}
\usepackage{amsfonts,amssymb,amsmath}
\usepackage{psfrag}
\usepackage[belowskip=2pt]{caption}
\usepackage{booktabs}
\usepackage{mdframed}

\biboptions{sort&compress}
\usepackage{tikz}                              
\usetikzlibrary{matrix,chains,positioning,decorations.pathreplacing,arrows}

\usepackage{dsfont}

\journal{Journal}

\begin{document}

\begin{frontmatter}

\title{Physics-based Machine Learning for Computational Fracture Mechanics}

\author{Fadi Aldakheel\corref{cor}\(^{a}\)}
\ead{fadi.aldakheel@ibnm.uni-hannover.de}
\ead[url]{https://www.ibnm.uni-hannover.de}
\cortext[cor]{Corresponding author.} 
\author{Elsayed S. Elsayed\(^{a}\)}
\author{Yousef Heider\(^{a}\)}
\author{Oliver Weeger\(^{b}\)}

\address{ \(^a\) Institute of Mechanics and Computational Mechanics, Leibniz Universit\"at Hannover, Appelstrasse 9a, 30167 Hannover, Germany} 

\address{ \(^b\) Cyber-Physical Simulation, Department of Mechanical Engineering, Technical University of Darmstadt, 64293 Darmstadt, Germany} 

\begin{abstract}
This study introduces a physics-based machine learning ($\phi$ML) framework for modeling both brittle and ductile fractures. Unlike {\it physics-informed neural networks}, which solve partial differential equations by embedding physical laws as soft constraints in loss functions and enforcing boundary conditions via collocation points, our framework integrates physical principles, such as the governing equations and constraints, directly into the neural network architecture. This approach eliminates the dependency on problem-specific retraining for new boundary value problems, ensuring adaptability and consistency. By embedding constitutive behavior into the network's foundational design, our method represents a significant step toward unifying material modeling with machine learning for computational fracture mechanics. Specifically, a feedforward neural network is designed to embed physical laws within its architecture, ensuring thermodynamic consistency.
Building on this foundation, synthetic datasets generated from finite element-based phase-field simulations are employed to train the proposed framework, focusing on capturing the homogeneous responses of brittle and ductile fractures. Detailed analyses are performed on the stored elastic energy and the dissipated work due to plasticity and fracture, demonstrating the capability of the framework to predict essential fracture features.
%
The proposed $\phi$ML framework overcomes the shortcomings of classical machine learning models, which rely heavily on large datasets (data-hungry) and lack guarantees of physical principles. By leveraging its physics-integrated design, the $\phi$ML framework demonstrates exceptional performance in predicting key properties of brittle and ductile fractures with limited training data. This ensures reliability, efficiency, and physical consistency, establishing a foundational approach for integrating machine learning with computational fracture mechanics.
\end{abstract}


\end{frontmatter}

\section{Introduction}
\label{sec1}

In response to the urgent societal and economic need for more sustainable and reliable mechanical systems with early fracture detection, computational mechanics is at a critical juncture. Accurately predicting the lifespan and health of structures under inelastic conditions is essential for ensuring the integrity and safety of critical structures, such as aircraft, buildings, or wind turbines. However, modeling rate-dependent inelastic material behaviors like plasticity, damage, or fracture entails significant challenges related to the complexity of the underlying physical processes and the limitations of current computational approaches:
\begin{itemize}
\item Time dependency: Accurately modeling path- or history-dependent material responses under varying loading conditions requires tackling numerical challenges related to selecting appropriate time integration schemes and formulating complicated material models.
\item Complex material behavior: This includes nonlinear and interdependent relationships between loads, stress, strain, and time.
\item Data scarcity: While experimental data is often limited and expensive to obtain, generating high-quality numerical data to address these challenges is also computationally demanding. This dual limitation restricts the development of accurate and reliable predictive models for fracture mechanics.
\end{itemize}
Conventional methods for predicting brittle and ductile fractures often depend on empirical models, which are constrained by assumptions on their structure or the need for extensive datasets, making them ineffective in scenarios with limited data or complex loading conditions.
Numerical techniques, including finite element and virtual element methods, have significantly advanced the modeling of fractures. Among these, the phase-field approach stands out as a powerful and versatile tool for simulating fracture processes. It offers a robust framework to capture crack initiation, propagation, and coalescence through a variational and thermodynamically consistent formulation. This makes it particularly well-suited for addressing both brittle and ductile fracture scenarios with high accuracy and reliability \cite{golahmar2023phase,abubakar2022influence,noii2022probabilistic,carrara2020framework,seiler2020efficient,storm2021comparative,baktheer2024phase,noii2022bayesian,chen2020adaptive,heider2020phase, Heider2021_review_PFHyd}.

Despite its strengths, the phase-field approach, like other computational fracture mechanics methods, is constrained by high computational demands.

\subsection{Machine learning methods}
As highlighted above, traditional numerical methods for fracture mechanics are struggling to keep up with the growing complexity of modern engineering systems, both mathematically and computationally. In this context, machine learning (ML) has emerged as a promising tool for advancing fracture mechanics. By leveraging its ability to recognize and learn intricate patterns from data, ML offers new opportunities for improving the efficiency and accuracy of fracture modeling.

Prior studies have explored the use of ML approaches for fracture prediction, including regression models, support vector machines (SVM), or random forests \cite{samaniego2020energy,aykol2021perspective,zhan2021machine,feng2021machine,liu2020machine,muller2021machine,zhan2021data}. Also deep learning has been widely applied to fracture mechanics, providing efficient and accurate solutions. Feed-forward neural networks (FFNNs) have been effectively utilized to model brittle fracture behavior, delivering results comparable to exact solutions under various loading conditions \cite{aldakheel2021feed}. Convolutional neural networks (CNNs) have demonstrated their capability for rapid corrosion prediction in coated magnesium alloys, significantly reducing computational costs while maintaining high predictive accuracy \cite{ma2025rapid}. Deep neural networks (DNNs) and deep reinforcement learning (DRL) have further advanced multi-scale computational fracture mechanics, offering flexible, data-driven approaches to model complex dependencies \cite{heider2021multi,ChaabanEtAl2023_ML_LBM}. Unsupervised learning, such as $k$-means clustering, has accelerated multiscale simulations of quasi-brittle materials by grouping integration points and reducing redundant computations \cite{chaouch2024unsupervised}.

While these ``black-box'' ML methods and data-driven approaches are often more efficient than classical numerical methods, their accuracy depends largely on the amount of available training data, i.e., they are \textbf{data hungry}. 
This can make them unreliable, unrobust, and uninterpretable, i.e., their accuracy cannot be guaranteed, they are sensitive to changes in the data they are trained on, and it is difficult to understand how they make predictions. Thus, they do not generalize well when evaluated outside the range of their training data \cite{montans2023}. All of these aspects are critical in fracture mechanics problems, where, on the one hand, only little experimental data are available or are expensive to generate, and, on the other hand, high demands on the accuracy and reliability of models and simulation results are required due to safety and certification regulations.

\subsection{Physics-based machine learning} \label{sec1-PhysicsML}
Thus, to deliver \textbf{physically valid, trustworthy, and ideally also explainable predictions}, generalizable ML models should be domain-aware, i.e., they should (1) be \emph{based on} physical principles, (2) \emph{embed} physical requirements into their architectures, (3) have a well-defined scope of applicability and generality, (4) provide quantifiable error measures, and (5) be evaluated against benchmark datasets \cite{brodnik2023}. To tackle these challenges,  only in the last 4 to 5 years researchers {from mathematics,  computer, and engineering sciences} have begun to \textbf{systematically combine the well-known principles of physical modeling and established numerical simulation methods with ML}, i.e., started developing \textbf{\emph{physics-based} machine learning} ($\phi$ML) methods.
In the literature, such approaches are also referred to as ``scientific'', ``structure-preserving'' or ``hybrid'' machine learning \cite{baker2019,karniadakis2021piml,as2023mechanics,michopoulos2024scientific}.

Most prominently, \emph{physics-informed} neural networks (PINNs) incorporate residual terms of differential equations into the loss function of the NN \cite{raissi2019,di2024physics}. PINNs have recently been explored in fracture mechanics, particularly for phase-field \emph{brittle fracture}, offering a data-driven alternative to traditional numerical methods. \citet{goswami2020transfer} introduced the first PINN framework for phase-field fracture modeling, simplifying boundary condition enforcement and leveraging lower-order derivatives for faster training. Despite potential computational savings, the need for retraining at each load step presents inefficiencies, which the authors mitigated using transfer learning to partially retrain the network. \citet{manav2024phase} employed the deep Ritz method (DRM) within PINNs to model phase-field fractures, capturing crack nucleation, propagation, branching, and coalescence. The method demonstrated robustness and close agreement with FEM solutions, showcasing PINNs ability to handle complex fracture scenarios effectively. \citet{ghaffari2023deep} compared multiple PINN variants, including PINNs, VPINNs, and VE-PINNs, in modeling brittle fracture. They analyzed boundary condition imposition, computational costs, and sensitivity to network parameters, highlighting the method's dependency on hyperparameter tuning and the computational burden associated with large-scale simulations.

While these approaches show potential, limitations such as computational inefficiency, sensitivity to hyper-parameters, and difficulties in achieving physical consistency underscore the need for a {\it stronger integration} of physical laws into machine learning frameworks. This motivates the development of physics-based machine learning, to overcome these challenges. So far, $\phi$ML models that aim to incorporate thermodynamic consistency have been proposed for hyperelasticity \cite{le2015,gonzalez2020,fernandez2021,klein2022}, as well as inelastic behaviors such as visco-elasticity \cite{tac2023a,rosenkranz2024a} or elasto-plasticity \cite{masi2021,vlassis2021,masi2023,fuhg2023,vlassis2023}. However, fully \emph{physics-augmented} inelastic formulations are yet to be extended to brittle and ductile fracture mechanics.

\medskip

Building on these foundations, this study presents a novel $\phi$ML framework tailored for modeling brittle and ductile fractures. By embedding governing equations and physical constraints directly into the neural network architecture, the framework captures the essential physics of fracture processes while maintaining computational efficiency. This design of the $\phi$ML framework ensures that it inherently respects physical laws, allowing it to deliver consistent and reliable predictions without requiring retraining for new boundary value problems or changing conditions. The framework is trained and validated using synthetic datasets derived from finite element-based phase-field simulations, specifically designed to represent homogeneous brittle and ductile fracture behaviors \cite{aldakheel2021simulation}. In this way, key fracture characteristics, such as energy degradation, plasticity and phase-field evolution, stress-strain response, and total dissipation, will be accurately predicted, showcasing the model's ability to replicate complex fracture phenomena. 

This paper is structured as follows: Section \ref{sec2} summarizes the governing equations in computational fracture mechanics, providing the mathematical basis for modeling material failure. Section \ref{sec3} provides numerical results of coupling elastoplasticity to fracture. These numerical data are needed for training both classical and physics-based ML models. Section \ref{sec4} starts by reviewing classical ML models, highlighting their applications and limitations in predicting fracture and inelastic behavior. Thereafter, a novel physics-based ML model is introduced, detailing how physical principles are embedded to enhance predictive accuracy and reliability in fracture mechanics problems. Finally, Section \ref{sec5} presents the results and a discussion, showcasing the model\'s performance in simulating brittle and ductile fracture scenarios, and demonstrating its potential for future applications in complex failure analysis and lifetime predictions. Finally, section \ref{sec6} presents concluding remarks and future aspects.

\section{Governing equations for computational fracture mechanics}
\label{sec2}
%
\begin{figure}[b]%
\centering
\includegraphics*[width=0.56\textwidth]{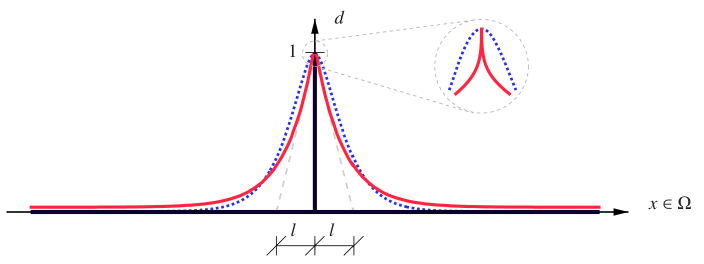} \qquad
\includegraphics*[width=0.34\textwidth]{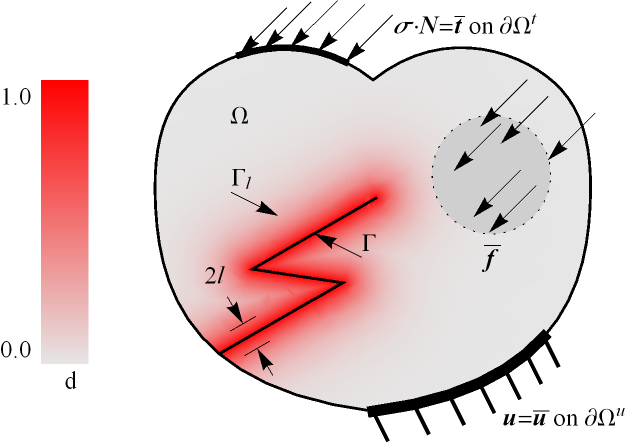}
\caption{
Left: One-dimensional cracked bar. The black thick line represents a sharp crack at $x=0$. The red and blue lines represent diffusive crack modeling. For the red line $\gamma_l=
d^2/2l + l\vert \nabla d \vert^2/2$ with regularization profile
$exp[ - \vert x \vert /l ]$ satisfying $d(0)=1$ and for the blue dotted line
$\gamma_l = d^2/2l + l |\nabla d|^2/4 + l^3(\Delta d)^2/32$ with
regularization profile $exp[ - 2\vert x \vert /l ] (1+ 2\vert x \vert
/l)$ satisfying $d(0)=1$ and $d^{\prime}(0)=0$. 
Right: Solid with a regularized crack and boundary conditions.}
\label{solid-FE}%
\end{figure}%

This section outlines a theory of phase-field fracture in elastic-plastic solids. It is based on a minimization of a pseudo-potential energy for the coupled problem. To this end, let $\Omega\subset{\calR}^{\delta}$ with $\delta = 1,2,3$ be a solid domain, as visualized in \rFig{solid-FE}. The response of a fracturing solid at material points $\Bx\in\Omega$ and time $t\in \calT = [0,T]$ is described by the displacement field $\Bu(\Bx,t)$ and the crack phase-field $d(\Bx,t)$ as
\begin{equation}
\Bu: 
\left\{
\begin{array}{ll}
 \Omega \times \calT \rightarrow \calR^{\delta} \\
 (\Bx, t)  \mapsto \Bu(\Bx,t)
\end{array}
\right.
\AND
d: 
\left\{
\begin{array}{ll}
 \Omega \times \calT \rightarrow [0,1] \\
 (\Bx, t)  \mapsto d(\Bx,t)
\end{array}
\right.
\WITH
\dot{d} \ge 0  \ .
\label{phi-d-fields}
\end{equation}
The crack phase-field $d(\Bx,t)=0$ and $d(\Bx,t)=1$ refer to the unbroken and fully broken state of the material, respectively, see \rFig{solid-FE}.

The gradient of the displacement field defines the symmetric strain tensor of the geometrically linear theory
\begin{equation}
\Bve = \nabla_s \Bu = \sym[ \nabla \Bu ] := \frac{1}{2} [\nabla\Bu + \nabla\Bu^T]
\ .
\end{equation}
The time-dependent Dirichlet- and Neumann boundary conditions of the solid are defined as
\begin{equation}
\Bu = \overline{\Bu} \ \textrm{on}\ \partial\Omega^u
\AND
\Bsigma \cdot \BN 
= \overline{\Bt} \ \textrm{on}\ \partial\Omega^{t}
\end{equation}
with a prescribed displacement and external traction on the surface $\partial \Omega = \partial \Omega^u\cup \partial
\Omega^{t}$ of the solid domain, where $\partial \Omega^{u} \cap \partial \Omega^{t}     = \emptyset$. The stress tensor $\Bsigma$ is the thermodynamic dual to $\Bve$.

To model irreversible, plastic deformations, the total strain tensor is additively decomposed into an elastic and a plastic part as
\begin{equation}
\Bve = \Bve^e + \Bve^p
\end{equation}
To account for the phenomenological hardening/softening response, we define the equivalent plastic strain variable by the evolution equation
\begin{equation}
\dot\alpha = \sqrt{\frac{2}{3}} \; \| \dot\Bve^p \|
\WITH
\dot\alpha \ge 0
\ 
\label{evol-alpha}
\end{equation}
as a local internal variable. The hardening variable starts to evolve from the initial condition $\alpha(\Bx,0) = 0$. For the phase-field problem, the sharp crack surface topology $\Gamma$ is regularized, yielding the crack surface functional $\Gamma_l$, defined as
\begin{equation}
\Gamma_l(d) = \int_{\Omega} \gamma_l(d, \nabla d) \, dV
\WITH
\gamma_l(d, \nabla d) =  
\dfrac{1}{2l} d^2 + \dfrac{l}{2} \vert \nabla d \vert^2 \,.
\label{gamma_l}
\end{equation}
As outlined in \cite{miehe+hofacker+schaenzel+aldakheel15}, this regularization is based on the definition of the crack surface density function $\gamma_l$ per unit volume of the solid and the fracture length scale parameter $l>0$ that governs the regularization, as illustrated in \rFig{solid-FE}. 
Evolution of \req{gamma_l} can be driven by the constitutive
functions as outlined in \cite{aldakheel16,aldakheel+wriggers+miehe17}, postulating a global evolution equation of regularized crack surface as
\begin{equation}
\vphantom{\int_{\Omega}}
\frac{d}{dt} \Gamma_l(d) =
\int_{\Omega} {\delta_d}
\gamma_l(d,\nabla d) \, \dot d \, dV := 
\frac{1}{l} \int_{\Omega} \big[ (1-d) \calH - \eta_d \dot{d} \big]\;
\dot{d} \, dV \ge 0 \ ,
\label{gamma-evol}
\end{equation}
where $\eta_d \ge 0$ is a material parameter that characterizes 
\textcolor{black}{the artificial/numerical viscosity of the crack propagation}. The variational derivative of $\gamma_l$ is defined as $\delta_d \gamma_l := \partial_d \gamma_l - \div [\partial_{\nabla d} \gamma_l] $. The crack driving force
$\calH$ is introduced as a {\it local history variable} that accounts on the irreversibility of the phase-field evolution, see \cite{miehe+aldakheel+raina16}. 
\newline
\newline
The above-introduced variables will characterize the brittle and ductile failure responses of a solid, based on the two  global primary fields
\begin{equation}
\mbox{Global Primary Fields}: \BfrakU := \{ \Bu, d \}
\label{global-fields}
\ .
\end{equation}
The constitutive approach to phase-field modeling focuses on the set
\begin{equation}
\mbox{Constitutive State Variables}: 
	\BfrakC := \{ \Bve^e, \alpha, d, \nabla d \}
\ .
\label{state}
\end{equation}
Following that, the work density function $W$ is split into an energetic elastic part and a dissipative plastic-fracture part as
\begin{equation}
{W}(\BfrakC) = 
{\Psi}_{e}(\Bve^e, d) + {D}{}_{pd} (\alpha, d, \nabla d) \ .
\label{pseudo-energy}
\end{equation}
This split assumes that the macroscopic elastic strain energy $\Psi_e$ is the only part of the total work density that is stored in the material. The constitutive expression for this part is
\begin{equation}
\Psi_e = (1-d)^2 \; \psi_e(\Bve^e) 
\qquad \mbox{with} \qquad
\psi_e(\Bve^e) = \frac{\kappa}{2}\; \mbox{tr}^2[\Bve^e] + \mu\;\tr [\dev(\Bve^e)^2]
\label{elas-part}
\end{equation}
in terms of the bulk modulus $\kappa > 0$ and the shear modulus $\mu > 0$. The accumulated dissipated work is defined as
\begin{equation}
{D}{}_{pd} (\alpha, d, \nabla d) = (1-d)^2 \; \big[ \psi_p(\alpha) - \psi_c \big] + \psi_c + 2 \frac{\psi_c}{\zeta} l~ {\gamma_l}(d, \nabla d)
\qquad \mbox{with} \qquad
\psi_p(\alpha) := y_0 \alpha + \frac{h}{2} \alpha^2\,.
\end{equation}
In this, $y_0$ is the initial yield stress, $h\ge 0$ is the isotropic hardening modulus, ${\psi}_c > 0$ is the critical fracture energy, and $\zeta$ is a fracture parameter that controls the post-critical range after crack initialization. These functions govern the constitutive expressions for the stress, the driving forces for plasticity and fracture
\begin{equation}
\Bsigma = \partial_{\Bve} {\Psi}{}_e \quad ,
\quad  \qquad
\Bf_p := - \partial_{\Bve^p} {\Psi}{}_e
\qquad \mbox{and} \qquad
f_d := - \partial_{d} {\Psi}{}_e \quad  ,
\label{stress-driving-forces}
\end{equation}
as well as the associated resistance functions for plasticity and fracture
\begin{equation}
r_p = \partial_{\alpha} {D}{}_{pd}
\qquad \mbox{and} \qquad
r_d = \partial_{d} {D}{}_{pd} \ .
\end{equation}
The evolution of the plastic strains and the fracture phase-field is governed by the two threshold functions
\begin{equation}
{\phi}{}_p(\Bf_p, r_p) = \big\| \dev[\Bf_p] \big\| - \hbox{$\sqrt{\frac{2}{3}}$}\;  r_p
\AND
{\phi}{}_d(f_d - r_d) = f_d - r_d
\ .
\end{equation}
Following \cite{aldakheel16}, we define the dissipation potential function for both plasticity and fracture dissipation mechanisms as
\begin{equation}
{V} (\dot{\BfrakC}) = 
\sup_{\Bf_p, r_p, f_d - r_d} \,  
\big[\; \Bf^p : \dot{\Bve}{}^p - r_p~ \dot{\alpha} + (f_d-r_d)~ \dot{d} 
- {V}^{\ast} (\Bf_p, r_p, f_d - r_d)
\; \big]
\ ,
\end{equation}
in terms of the {dual dissipation potential function}
\begin{equation}
{V}^{\ast} (\Bf_p, r_p, f_d - r_d) = 
\frac{1}{2\eta_p} \big\langle {\phi}{}_p(\Bf_p, r_p) \big\rangle^2 +
\frac{1}{2\eta_d} \big\langle {\phi}{}_d(f_d - r_d) \big\rangle^2 \ ,
\end{equation}
where $\langle x \rangle := ( x + \vert x\vert )/2$ is the Macaulay bracket.
$\eta_p$ and $\eta_d$ are additional material parameters that characterize the viscosity of the plastic deformation and the crack propagation. The necessary conditions of the optimization problem determine the {plastic flow rules} and the evolution equation for the crack phase-field are
\begin{equation}
\dot\Bve^p = \lambda_p \, \partial_{\Bf_p} {\phi}{}_p \quad ,
\quad \qquad
\dot\alpha = - \lambda_p \, \partial_{r_p} {\phi}{}_p
\qquad \mbox{and} \qquad
\dot d = \lambda_d \, \partial_{f_d-r_d} {\phi}{}_d \quad ,
\end{equation}
along with the two loading-unloading conditions 
\begin{equation}
\lambda_p:= \frac{1}{\eta_p} \left\langle {\phi}{}_p \right\rangle \ge 0 \;\, ; \;\ {\phi}{}_p \le 0 \;\, ; \;\ \lambda_p {\phi}{}_p = 0
\quad \AND \quad
\lambda_f:= \frac{1}{\eta_d} \left\langle {\phi}{}_d \right\rangle \ge 0 \;\, ; \;\ {\phi}{}_d \le 0 \;\, ; \;\ \lambda_f {\phi}{}_d = 0
\end{equation}
of the plastic and fracture response, respectively. 

Using the quantities introduced above, we define the two governing partial differential equations (PDEs) of the coupled problem:
\begin{mdframed}
\begin{itemize}
    \item The solid $\Omega$ has to satisfy the equation of equilibrium as
\begin{equation}
{
\div\,\Bsigma + \overline\Bf = \Bzero\,,
}
\end{equation}
where dynamic effects are neglected, $\Bsigma$ is the stress tensor and $\overline\Bf$ is the given body force.
\item The evolution statement
\req{gamma-evol} provides the local equation for the
evolution of the crack phase-field in the domain $\Omega$ along with its homogeneous Neumann boundary condition as
\begin{equation}
\begin{aligned}
\eta_d \dot{d}
&= 
(1-d) {\calH} 
-
[ \, d - l^2 \Delta d \, ] \quad & \mbox{in}~\Omega\times\calT \ ,
\\
\nabla d \cdot \BN &= 0
\quad & \mbox{on} ~ \partial\Omega\times\calT \ .
\end{aligned}
\label{phase-field}
\end{equation}
The history field $\calH$ is defined in line with \cite{miehe+aldakheel+raina16,miehe+kienle+aldakheel+teichtmeister16} by
\begin{equation}
\calH := \max_{s\in [0,t]} \calD(\Bve^e,\alpha ; s) \ge 0
\WITH
\calD := \zeta \Big< {\psi}_e + {\psi}_{p} - \psi_c \Big>_+\,.
\label{history-field}
\end{equation}
Herein, the Macaulay brackets $\langle x \rangle_+ := (x + \vert x\vert)/2$ ensure the irreversibility of the crack evolution.
\end{itemize}
\end{mdframed}

\begin{figure}[b]%
\centering
\begin{subfigure}{0.49\linewidth}
\includegraphics*[width=\textwidth,trim={10mm 5mm 18mm 18mm},clip]{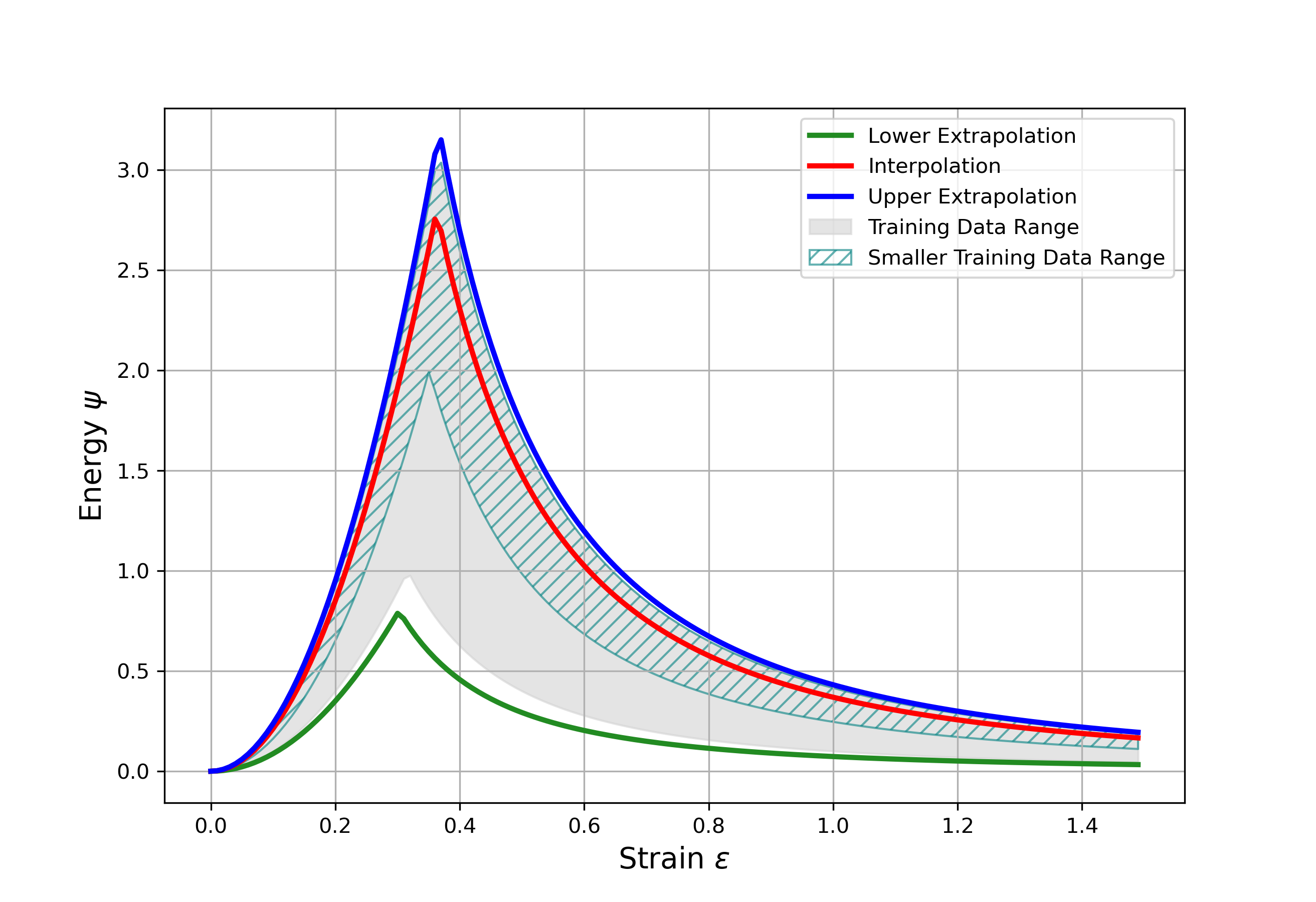} 
\caption{Energy}
\end{subfigure}
\hfill
\begin{subfigure}{0.49\linewidth}
\includegraphics*[width=\textwidth,trim={10mm 5mm 18mm 18mm},clip]{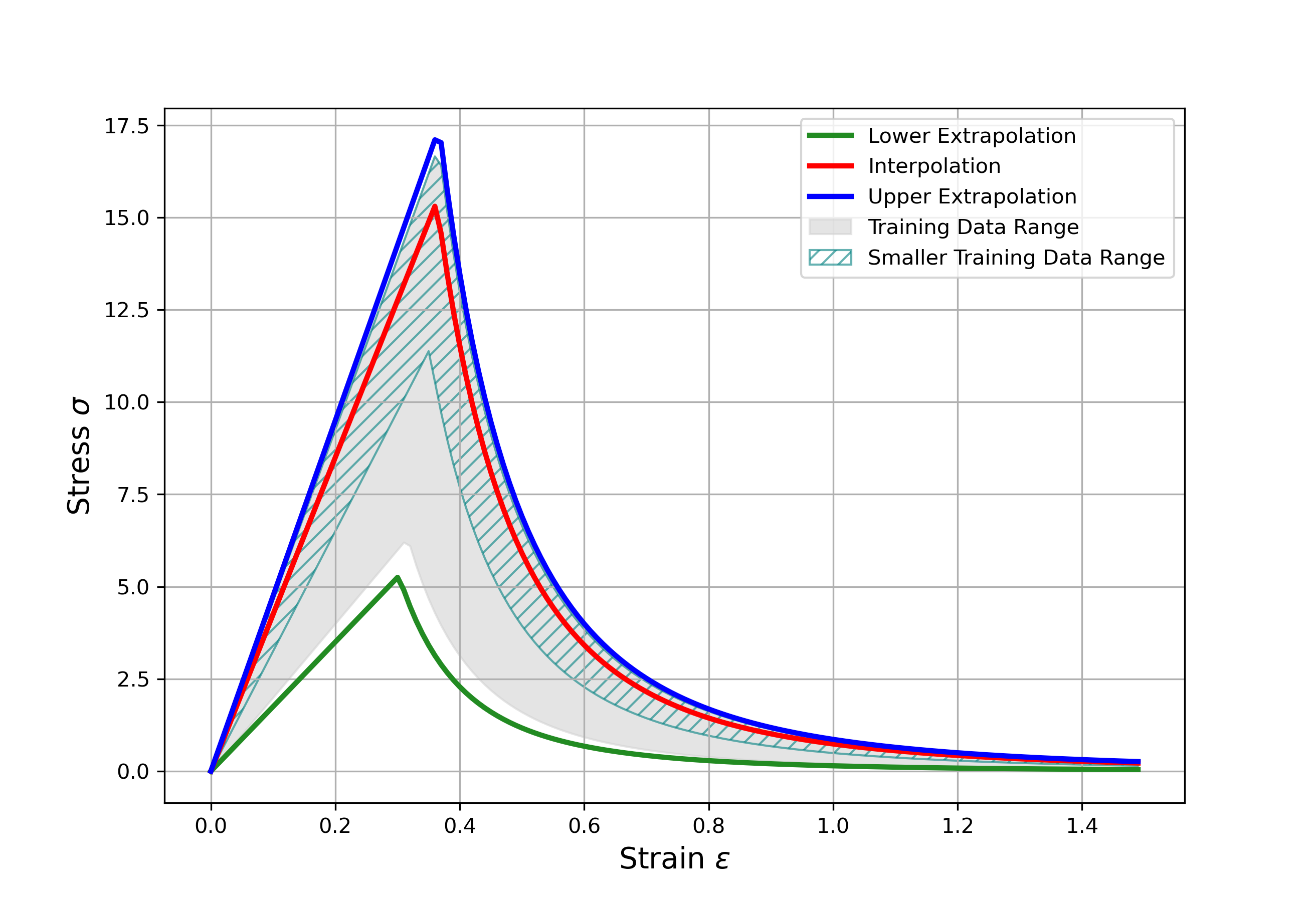}
\caption{Stress}
\end{subfigure}
\\[1mm]
\begin{subfigure}{0.49\linewidth}
\includegraphics*[width=\textwidth,trim={10mm 5mm 18mm 18mm},clip]{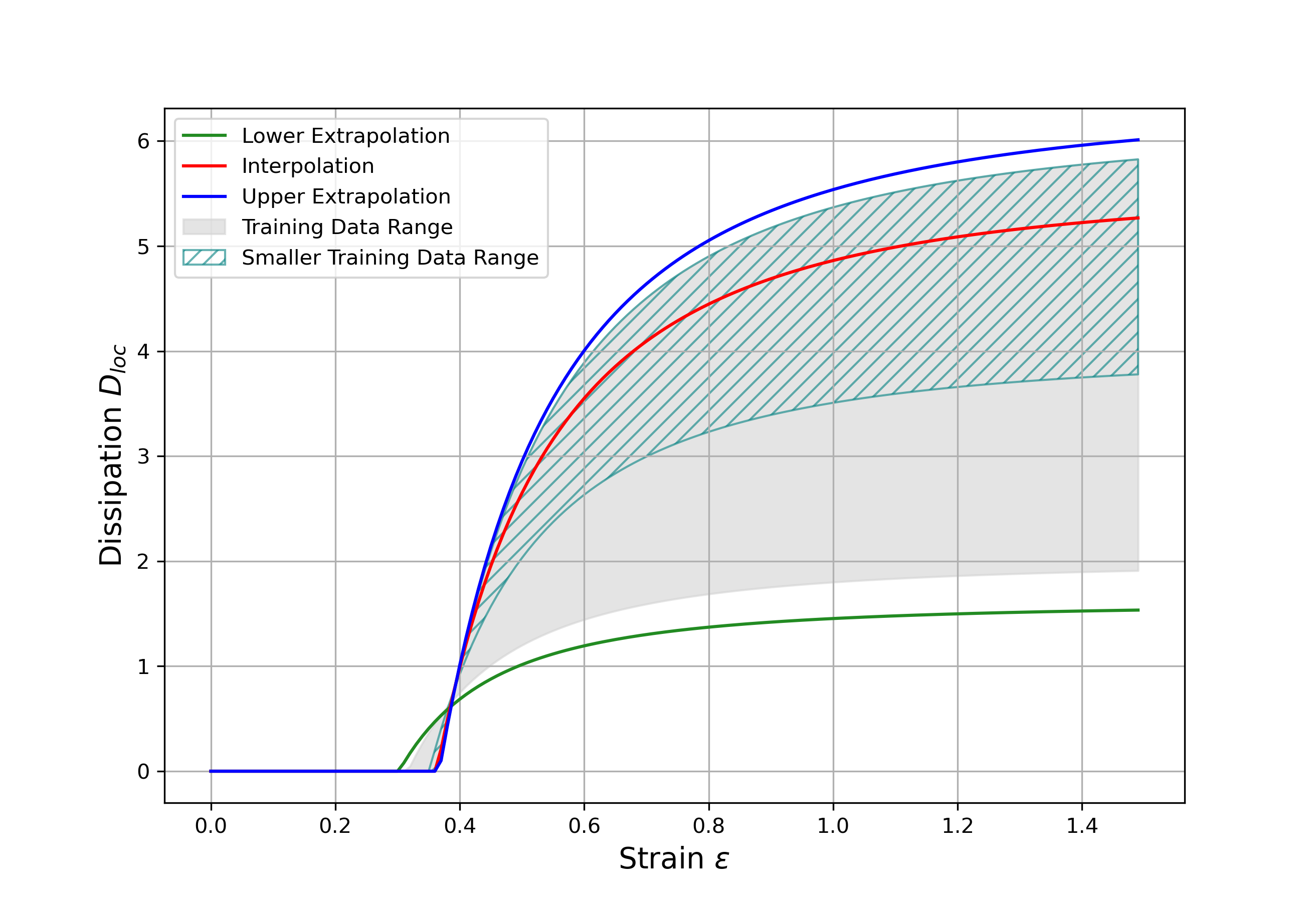} 
\caption{Dissipation}
\end{subfigure}
\hfill
\begin{subfigure}{0.49\linewidth}
\includegraphics*[width=\textwidth,trim={10mm 5mm 18mm 18mm},clip]{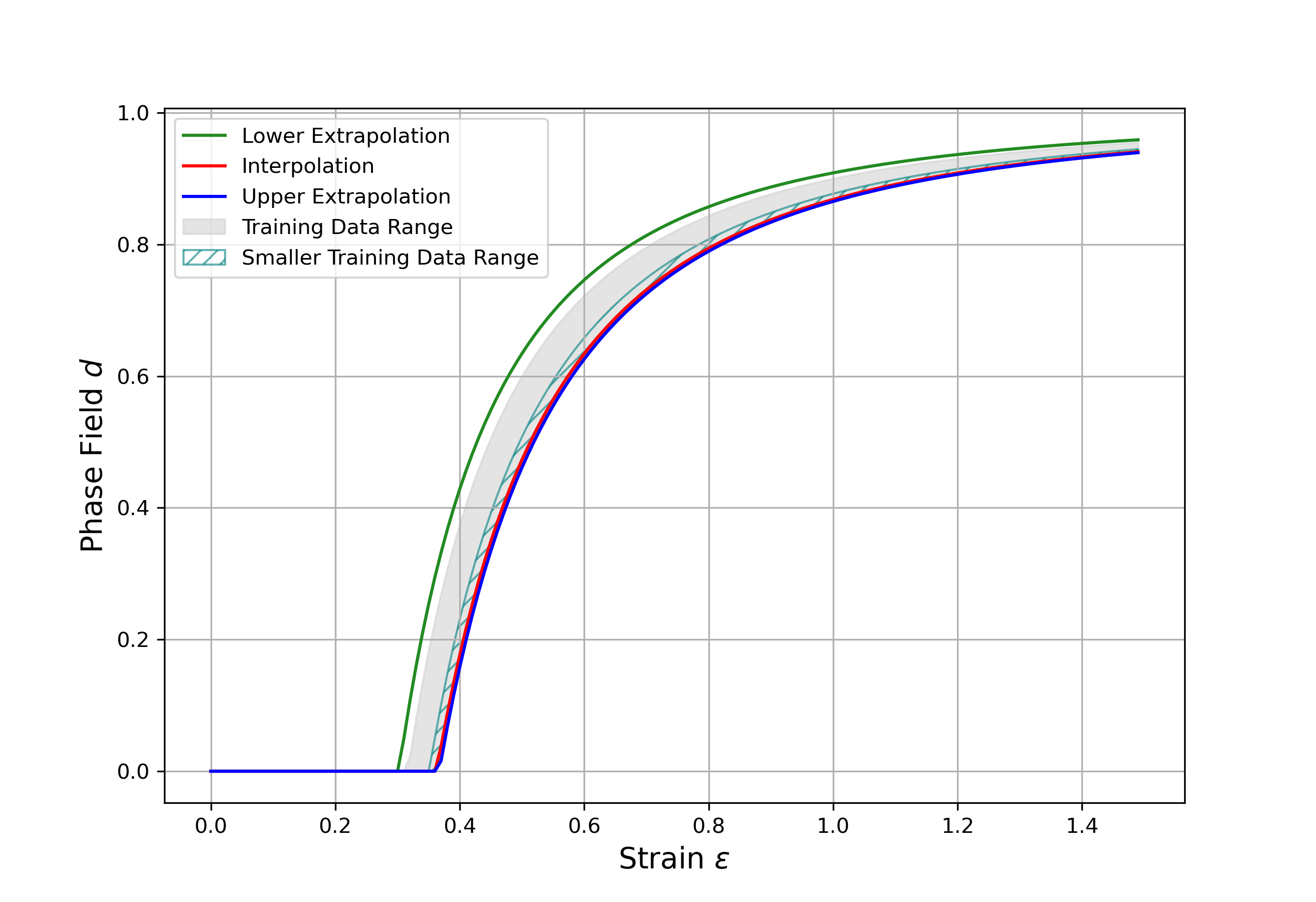}
\caption{Phase-field}
\end{subfigure}
\vspace*{-2mm}
\caption{
Brittle fracture: Data generation based on different values of Young's modulus and critical fracture energy.}
\label{data-generationbrit-brittle}%
\end{figure}%

\begin{figure}[t]%
\centering
\begin{subfigure}{0.49\linewidth}
\includegraphics*[width=\textwidth,trim={9mm 5mm 18mm 18mm},clip]{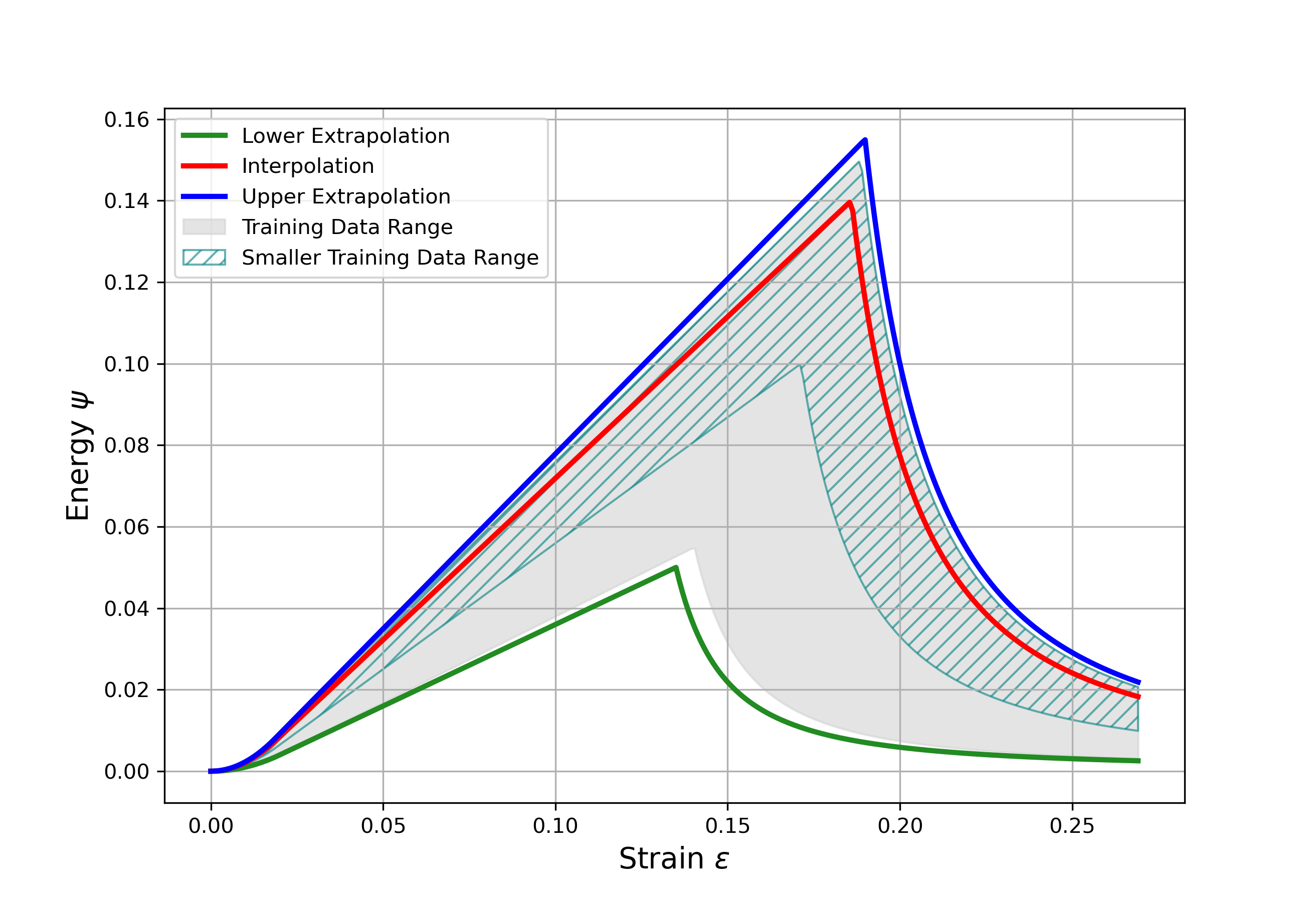}
\caption{Energy}
\end{subfigure}
\hfill
\begin{subfigure}{0.49\linewidth} 
\includegraphics*[width=\textwidth,trim={9mm 5mm 18mm 18mm},clip]{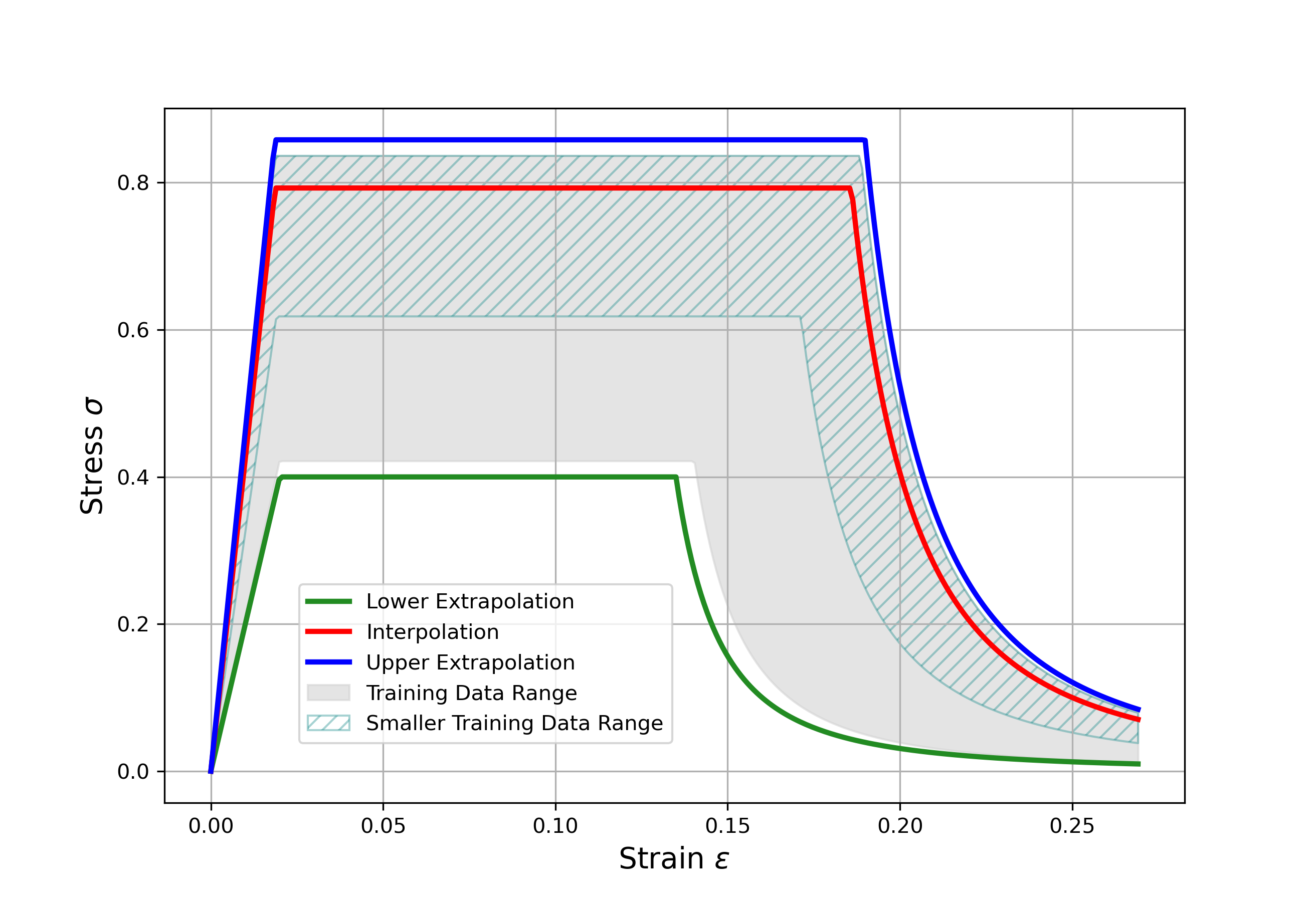}
\caption{Stress}
\end{subfigure}
\\[1mm]
\begin{subfigure}{0.49\linewidth} 
\includegraphics*[width=\textwidth,trim={9mm 5mm 18mm 18mm},clip]{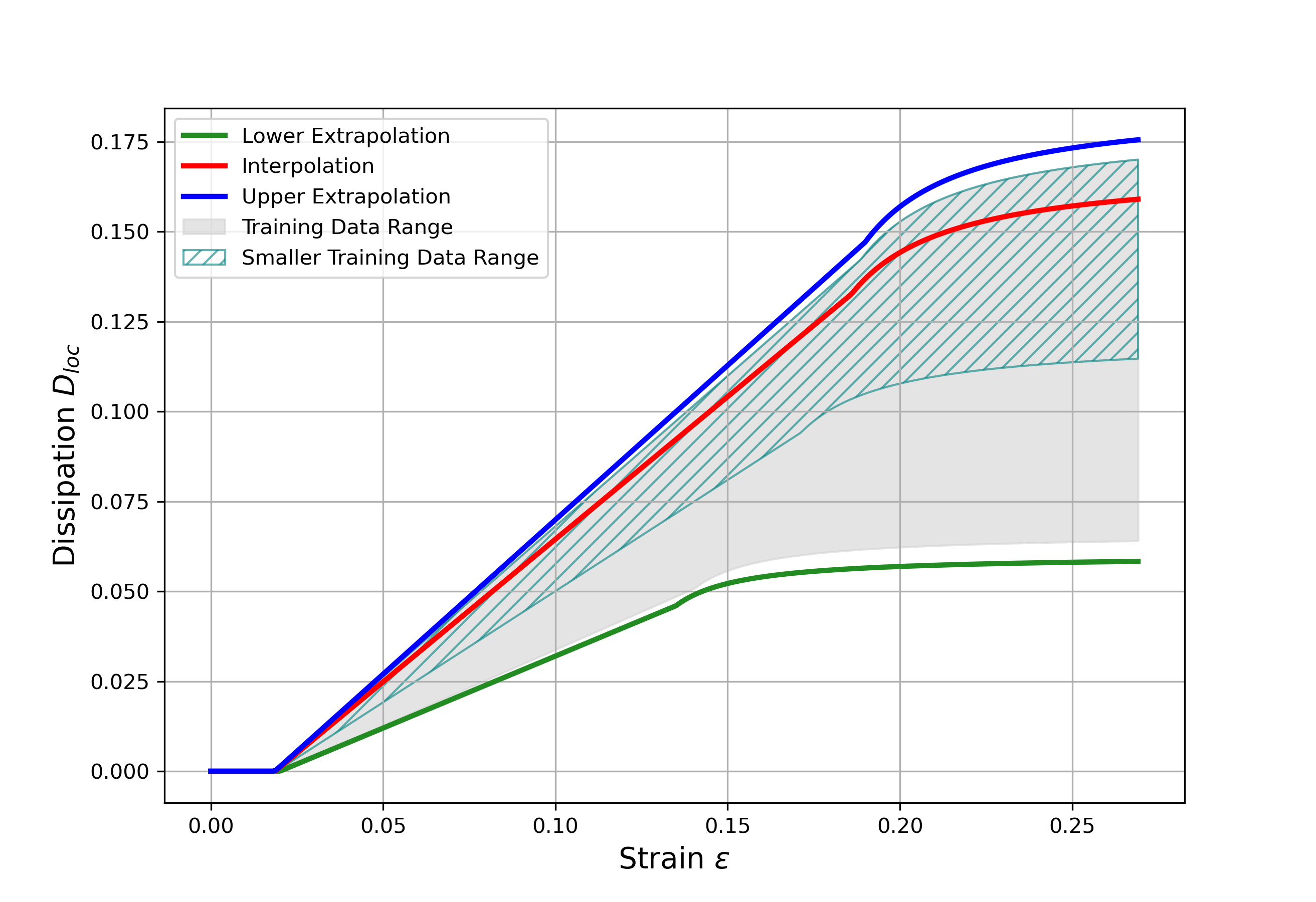}
\caption{Dissipation}
\end{subfigure}
\hfill
\begin{subfigure}{0.49\linewidth}
\includegraphics*[width=\textwidth,trim={9mm 5mm 18mm 18mm},clip]{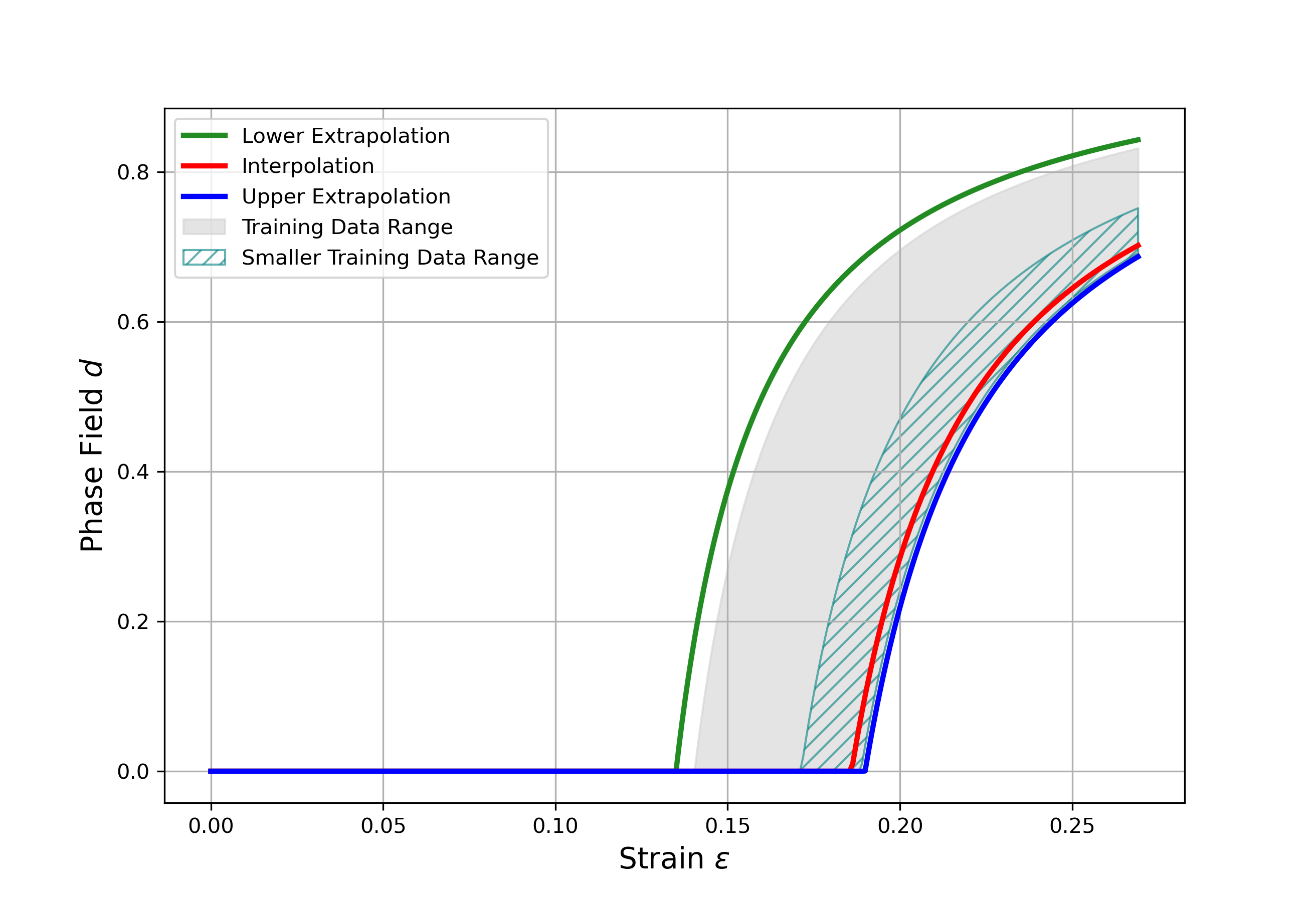}
\caption{Phase-field}
\end{subfigure}
\vspace*{-2mm}
\caption{
Ductile fracture: Data generation based on different values of Young's modulus, yield stress, and critical fracture energy.}
\label{data-generationbrit-ductile}%
\end{figure}%
%

\section{Model investigations and data generation}
\label{sec3}

The current work is the first to explore combining physics-based and ML models to address brittle and ductile fracture mechanics problems. Consequently, the focus is on a simplified, homogeneous response for both brittle (elastic-fracture) and ductile (elastic-plastic-fracture) solids. 
This section presents numerical results on the coupling of elastoplasticity and fracture, which will subsequently serve as training data for both naive ML (purely data-driven) and $\phi$ML models.

The material parameters used in the simulations are the same as in the reference work of \cite{miehe+aldakheel+raina16}. For data generation required for ML models, the elastic Young's modulus $E$ varies in the range $\{20 - 50\}$~GPa, the plastic yield stress $y_0$ (bounds the effective stress) is in the range $\{0.4 - 0.85\}$~GPa, the critical fracture energy $\psi_c$ (determines the onset of fracture) is in the range $\{0.05 - 0.155\}$~GPa, and the fracture parameter that controls the shape of the softening due to fracture is set to $\zeta=1$. In this work, ideal plasticity theory is considered, i.e., the isotropic hardening modulus is set to $h=0$. 

The brittle and ductile fracture model responses are visualized as gray-shaded regions in Figures \ref{data-generationbrit-brittle} and \ref{data-generationbrit-ductile}, respectively. Hereby, the elastic strain energy, defined in \req{elas-part}, is demonstrated in Figure \ref{data-generationbrit-brittle}a) for brittle fracture and Figure \ref{data-generationbrit-ductile}a) for ductile fracture. The stresses defined in \req{stress-driving-forces} are shown in Figure \ref{data-generationbrit-brittle}b) for brittle fracture and in Figure \ref{data-generationbrit-ductile}b) for ductile fracture, whereas the phase-field evolution defined in \req{phase-field} is shown in Figure \ref{data-generationbrit-brittle}d) for brittle fracture and in Figure \ref{data-generationbrit-ductile}d) for ductile fracture. 
The accumulated dissipative work $D$ is decomposed into a fracture contribution $D_d$ for brittle materials and into both a plasticity contribution $D_p$ and a fracture contribution $D_d$ for ductile materials, as shown in Figures \ref{data-generationbrit-brittle}c) and \ref{data-generationbrit-ductile}c), respectively. This decomposition is defined as 
\begin{equation}
D= D_p  + D_d 
\qquad
\mbox{with}
\quad
D_p := \int_0^{\varepsilon_p} f_p  \; d \tilde\varepsilon_p 
\quad \AND \quad
D_d := \int_0^d f_d  \; d \tilde{d} \quad ,
\end{equation}
which are obtained by numerical integration.
\begin{figure}[t]%
\centering
\includegraphics*[width= 0.35\textwidth]{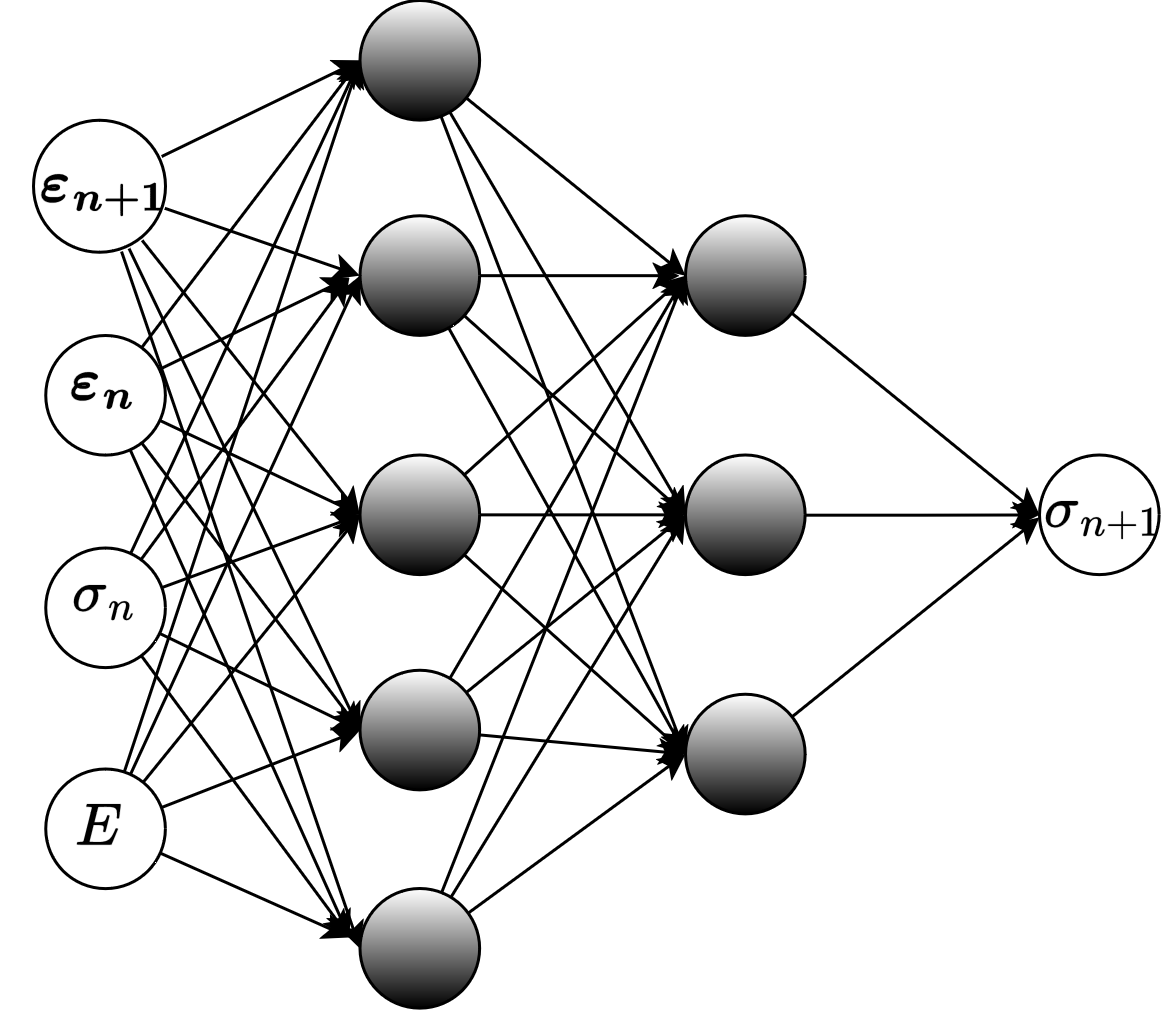} \;
\caption{Model architecture of \emph{uniformed} FFNN for fracture mechanics.}
\label{Classical-ML}%
\end{figure}%
%

\section{Physics-based machine learning model}
\label{sec4}
Future engineering progress depends on solving problems that require accurate and efficient computational modeling and simulations of failure mechanisms. To mitigate the high computational cost of solving fracture mechanics problems, as discussed in Section \ref{sec2}, ML is employed as an efficient approach to replace direct numerical calculations. 

\paragraph*{\underline{\textbf{Naive FFNN}}}
As a naive, purely data-driven ML approach, a standard Feed-Forward Neural Network (FFNN) is utilized to predict specific outputs based on given inputs. This study focuses on predicting stress responses using multi-channel inputs, including current and previous strain values, previous stress state, and Young's modulus. As illustrated in Figure \ref{Classical-ML}, this classical ML-based approach is applied to solve the brittle and ductile fracture model problems introduced in Section 
\ref{sec2}. This classical FFNN depicted in Figure \ref{Classical-ML} consists of two fully connected hidden layers, containing 16 and 8 neurons respectively with ReLU as an activation function, mapping the four input features to the current predicted stress state: 
\begin{equation}
    FFNN^{\sigma} : \quad  
    (\varepsilon_{n+1}, \varepsilon_n, \sigma_n, E)  
    \quad \longmapsto \quad  
    \sigma^{n+1}.
\end{equation}
Training is performed using the MSE loss function $\mathcal{L}_\sigma := \text{MSE}(\sigma)$ and Adam optimization approach.

As is well known, applying classical ML is an empirical and highly iterative process that requires training several models to achieve a satisfactory performance. This involves systematically testing different hyperparameter combinations to optimize performance. While such \emph{uninformed} ML approaches may accelerate computations by replacing direct numerical simulations, they remain data-hungry and often lack interpretability (black-box models). Thus, despite offering a more efficient alternative, their reliance on large datasets and extensive training presents new challenges, requiring careful consideration to balance accuracy, efficiency, and reliability in engineering applications. 

\paragraph*{\underline{\textbf{$\phi$ML model}}}
To overcome these limitations, ML models must integrate physical principles to enhance reliability and interpretability, as discussed in Section \ref{sec1-PhysicsML}. In this contribution, we argue that the reliable use of ML in fracture mechanics requires physics-based ML approaches that are guided by physical model information to preserve its essential structures and requirements. This ensures reasonable, reliable, and robust predictions. Furthermore, $\phi$ML models require manageable computational effort and training resources and are designed to operate effectively even with limited data availability. 

\begin{figure}[t]%
\centering
\includegraphics*[width=0.9 \textwidth]{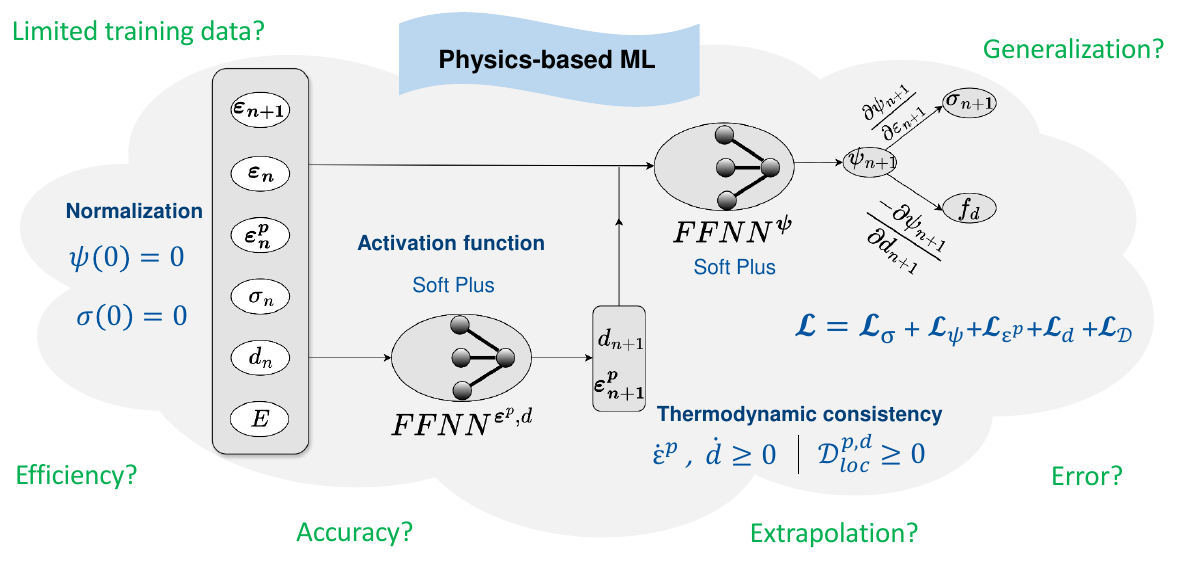} \;
\caption{Model architecture for \emph{physics-based} machine learning approach to fracture mechanics.}
\label{Physics-ML}%
\end{figure}

Thus, we develop a $\phi$ML architecture for fracture mechanics, which is presented in Figure \ref{Physics-ML}. This extended FFNN structure includes physical laws within the model architecture, ensuring normalization and thermodynamic consistency. The proposed architecture aligns with the recent works of \cite{masi2021thermodynamics,rosenkranz2023}, which integrate physical laws into ML models for elasto-plasticity. While these efforts represent significant progress, our approach extends beyond elastoplasticity to brittle and ductile fracture responses, where the interaction between elastoplastic deformation and crack evolution introduces additional complexities. By capturing these coupled mechanisms, our framework broadens the applicability of physics-based ML in failure modeling. Herein, the $\phi$ML model includes two separate, yet dependent neural networks that are implemented together: 
\begin{itemize}
    \item The first network $FFNN^{\varepsilon^p, d}$ receives as input the current and previous strain, previous stress, Young's modulus, and the previous values of both the plastic strain and phase-field. It then predicts the current plastic strain and phase-field as   
\begin{equation}
    FFNN^{\varepsilon^p, d} :\quad  
    (\varepsilon_{n+1}, \varepsilon_n, \varepsilon^p_n, \sigma_n, d_n, E)  
    \quad \longmapsto \quad  
    (\varepsilon^p_{n+1}, d_{n+1})
\end{equation}
\item The second network $FFNN^\psi$ takes the predicted phase-field and plastic strain from the first network together with the original inputs to compute the free energy, as follows:
\begin{equation}
    FFNN^{\psi} :\quad  
    (\varepsilon_{n+1}, \varepsilon_n, \varepsilon^p_n, \varepsilon^p_{n+1}, \sigma_n, d_n, d_{n+1}, E)  
    \quad \longmapsto \quad  
    \psi^{n+1}.
\end{equation}
\end{itemize}
The training process is formulated as a single optimization problem for both networks simultaneously, ensuring a unified learning framework. This approach is consistent with the \emph{thermodynamics informed} methodology of \cite{masi2021thermodynamics}, leading to the following loss function:
\begin{equation}
    \mathcal{L}_{total}:= \mathcal{L}_{\sigma} + \mathcal{L}_{\psi} + \mathcal{L}_{\varepsilon^p} + \mathcal{L}_d + \mathcal{L}_{D} \ ,
\end{equation}
where each term is defined as
\begin{equation}
    \mathcal{L}^\sigma := \text{MSE}(\sigma) \qquad \text{with} \qquad \sigma = \frac{\partial \psi_{n+1}}{\partial \varepsilon_{n+1}}
\end{equation}

\begin{equation}
    \mathcal{L}_\psi := \text{MSE}(\psi) \ , \qquad \mathcal{L}_{d} := \text{MSE}(d) \ , \qquad  \mathcal{L}_{\varepsilon^p} := \text{MSE}(\varepsilon^p) \ , \qquad \mathcal{L}_{{D}} := \text{MSE}({D}) \ .
\end{equation}
To enforce physically consistent dissipation, the total dissipation term is expressed as:
\begin{equation}
    {D} = \text{RELU}({D_d}) + \text{RELU}({D_p})\ ,
\end{equation}
where the fracture and plastic dissipation components evolve according to:
\begin{equation}
     D_{d}^{n+1} = D_{d}^{n} + \left(-\frac{\partial \psi_{n+1}}{\partial d_{n+1}}\right) \cdot (d_{n+1}-d_{n})
\quad \AND \quad
    D_{p}^{n+1} = D_{p}^{n} + \sigma_{n+1} \cdot (\varepsilon^p_{n+1} - \varepsilon^p_{n}) \ .
\end{equation}
In this work, all loss function terms are assigned equal weight during training, eliminating the need for a value representing the weight of each term in the total loss. Unlike the approaches in \cite{masi2021thermodynamics,rosenkranz2023}, which introduced separate weight terms for each loss component, in this work the same weight value of $1$ is assigned to all the loss terms on the loss function, allowing an equal contribution from all terms. Both subnetworks in our model employ two hidden layers with 16 and 8 neurons, using the \text{SoftPlus} activation function. For the neuron predicting the phase-field, a modified \text{ReLU}$^d$ function is applied to enforce an upper bound of $1$, ensuring the physical consistency of the phase-field variable $d$. This modification is defined as 
\begin{equation}
\text{ReLU}^d(x):= \begin{cases} 
0, & x \leq 0 \\
x, & 0 < x < 1 \\
1, & x \geq 1
\end{cases}
\end{equation}

\begin{figure}[b]%
\centering
 \includegraphics*[width=0.9\textwidth]{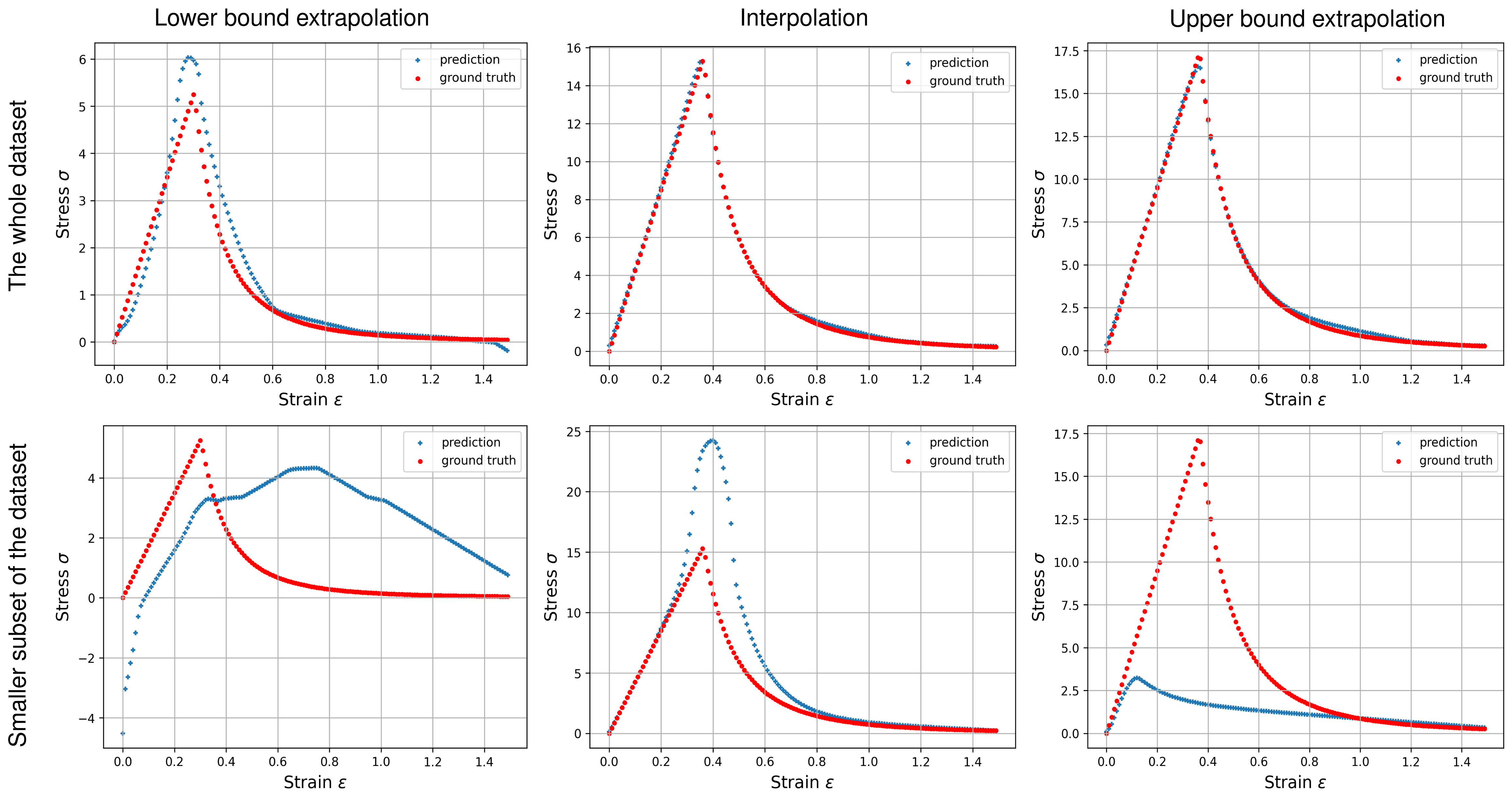} 
\caption{
Brittle fracture: Evaluation of purely data-driven, naive FFNN models using the full (top) versus reduced datasets (bottom) for the 3 test load paths (left: {\it lower bound extrapolation}, middle: {\it interpolation}, right column: {\it upper bound extrapolation}).
}
\label{brit-classical-ML}%
\end{figure}%

\begin{figure}[b]%
\centering
\includegraphics*[width=0.99\textwidth]{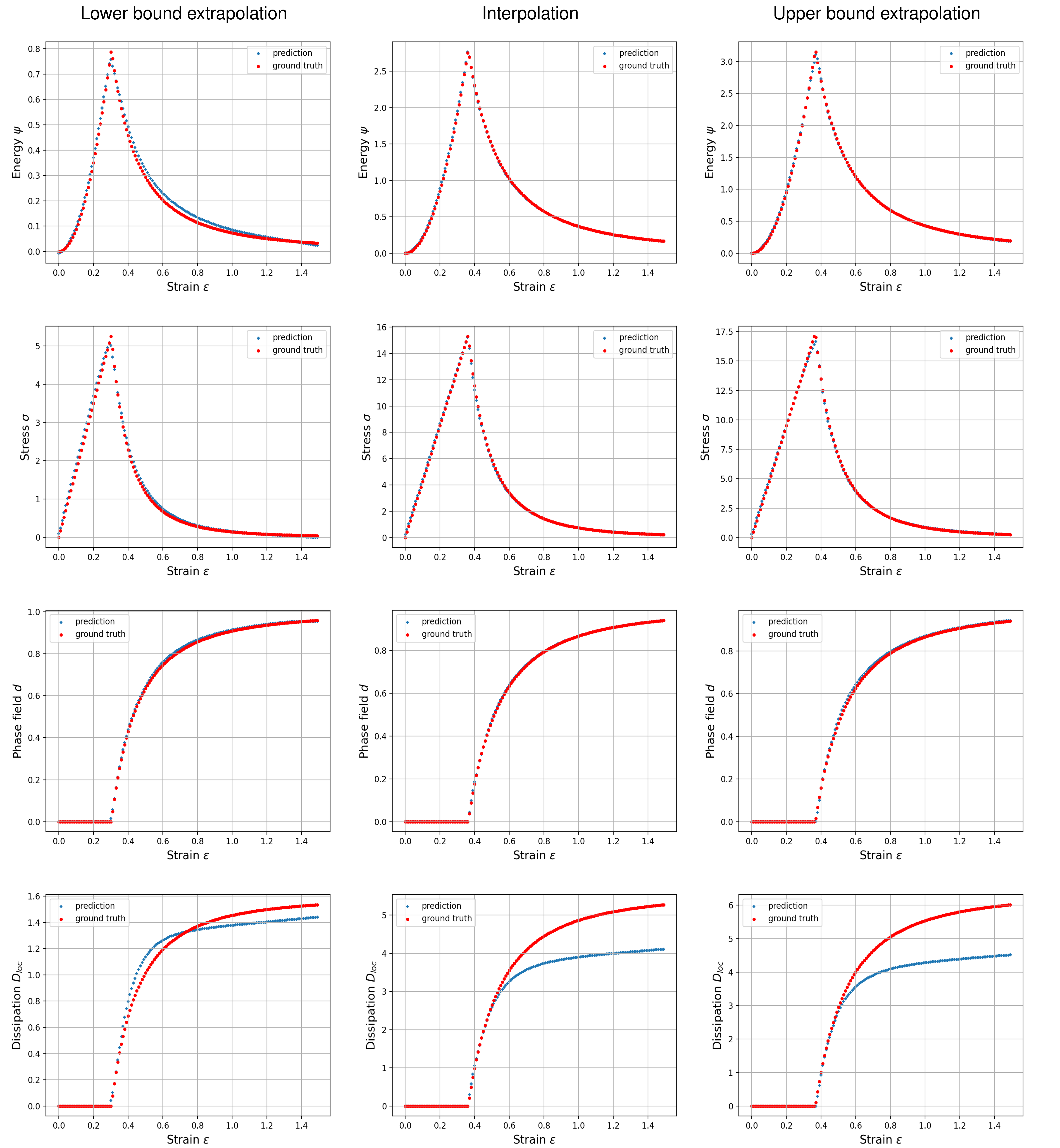} 

\caption{
Brittle fracture: Evaluation of $\phi$ML model using the full dataset for the 3 test load paths (left, middle, right column).
}
\label{brit-pML-full}%
\end{figure}%

\begin{figure}[b]%
\centering
 \includegraphics*[width=0.99\textwidth]{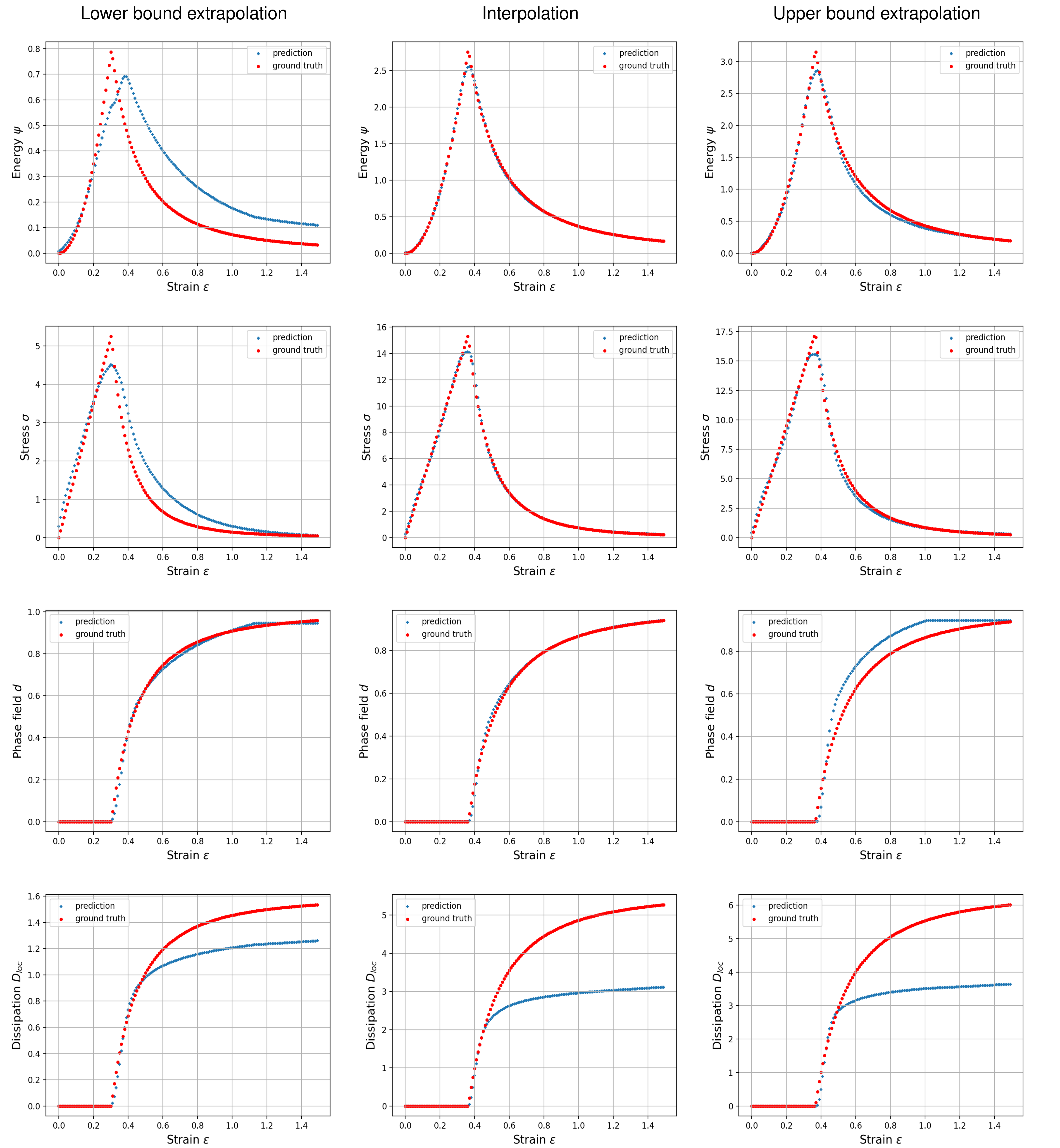} 

\caption{
Brittle fracture: Evaluation of $\phi$ML model using the reduced dataset for the 3 test load paths (left, middle, right column).
}
\label{brit-pML-reduced}%
\end{figure}%

\paragraph*{\underline{\textbf{Training and test data}}}
As outlined in Section \ref{sec3}, for dataset generation, the phase-field model is utilized with varying material parameters to train, validate, and test the ML models. The corresponding data are shown in Figure \ref{data-generationbrit-brittle} for brittle fracture and Figure \ref{data-generationbrit-ductile} for ductile fracture. To examine the effect of data availability on ML prediction accuracy, two dataset configurations are considered for model training, validation, and testing:
\begin{itemize}
\item \textbf{Full dataset}: Consisting of 23 cases, this dataset is used for training (18 load-path cases), validating (2 load-path cases), and testing (3 load-path cases) the ML model. The training dataset is represented by the gray-shaded regions in Figures \ref{data-generationbrit-brittle} and \ref{data-generationbrit-ductile}.

\item \textbf{Reduced dataset}: Represented by the hashed light-green area, this subset consists of 13 cases (half of the full training and validation dataset) selected to assess the model's performance with limited training data, see Figures \ref{data-generationbrit-brittle} and \ref{data-generationbrit-ductile}. {Specifically, 9 load-path cases are used for training, while 1 load-path case is allocated for validation and the same 3 load-path cases for testing as in the case of the full dataset.}

\end{itemize}
For both datasets, the ML models are {\it tested} under three scenarios (3 load-path cases) in Figures \ref{data-generationbrit-brittle} and \ref{data-generationbrit-ductile}, as follows: 
\begin{itemize}
\item Lower extrapolation (green line),
\item Interpolation (red line),
\item Upper extrapolation (blue line).
\end{itemize}

\section{Results and discussion}
\label{sec5}

In this section, the phase-field fracture numerical results (Section \ref{sec3}) will be used to train and evaluate both classical ML models and $\phi$ML models (Section \ref{sec4}). Furthermore, their predictive performance and generalization will thoroughly be compared. This comparison will highlight the benefits of integrating physical principles into ML for modeling future complex fracture scenarios.

\subsection{Brittle fracture}

For the brittle failure mechanism (E-F), the material remains elastic (E) up to a certain threshold $\psi_c$, after which it fractures (F), leading to final failure and separation. The numerical results from the E-F framework are used to generate the dataset for the ML models, as plotted in Figure \ref{data-generationbrit-brittle}. Hereby, a dataset of 3450 loading steps (23 load paths $\times$ 150 loading steps) is considered, to capture stress distribution, crack initiation, and propagation at each step.

Figure \ref{brit-classical-ML} illustrates the performance of purely data-driven, naive FFNN models for the full dataset (top row) and reduced dataset (bottom row), in predicting stress-strain relationships for brittle fracture under three {\it tested} dataset scenarios: lower bound extrapolation, interpolation, and upper bound extrapolation. The three ML test load-paths predictions are represented by blue dots, while the ground truth values are shown as red dots. While the full dataset model achieves reasonable accuracy during interpolation, it struggles in extrapolation scenarios, even with extensive data. The reduced subset further degrades performance across all scenarios, especially in extrapolation, highlighting the significant data dependency of data-driven ML models. These results demonstrate the inherent limitations of purely data-driven approaches, particularly their inability to generalize beyond the training range. They also emphasize the need for physics-informed ML models to enhance extrapolation capabilities and robustness.

Figures \ref{brit-pML-full} and \ref{brit-pML-reduced} demonstrate the strong performance of the physics-based ML model in predicting energy $\psi$, stress $\sigma$, phase-field $d$, and dissipation $D$ across lower bound extrapolation, interpolation, and upper bound extrapolation. For the full dataset (Figure \ref{brit-pML-full}), the model achieves near-perfect agreement with the ground truth for all tested dataset scenarios, including extrapolation. Even with a reduced subset (Figure \ref{brit-pML-reduced}), the model maintains accurate predictions, with minor deviations in extrapolation for stress and dissipation. These results highlight the robustness and reliability of the physics-based approach, outperforming the purely data-driven, naive FFNN model.


\begin{figure}[b]%
\centering
\includegraphics*[width=0.9\textwidth]{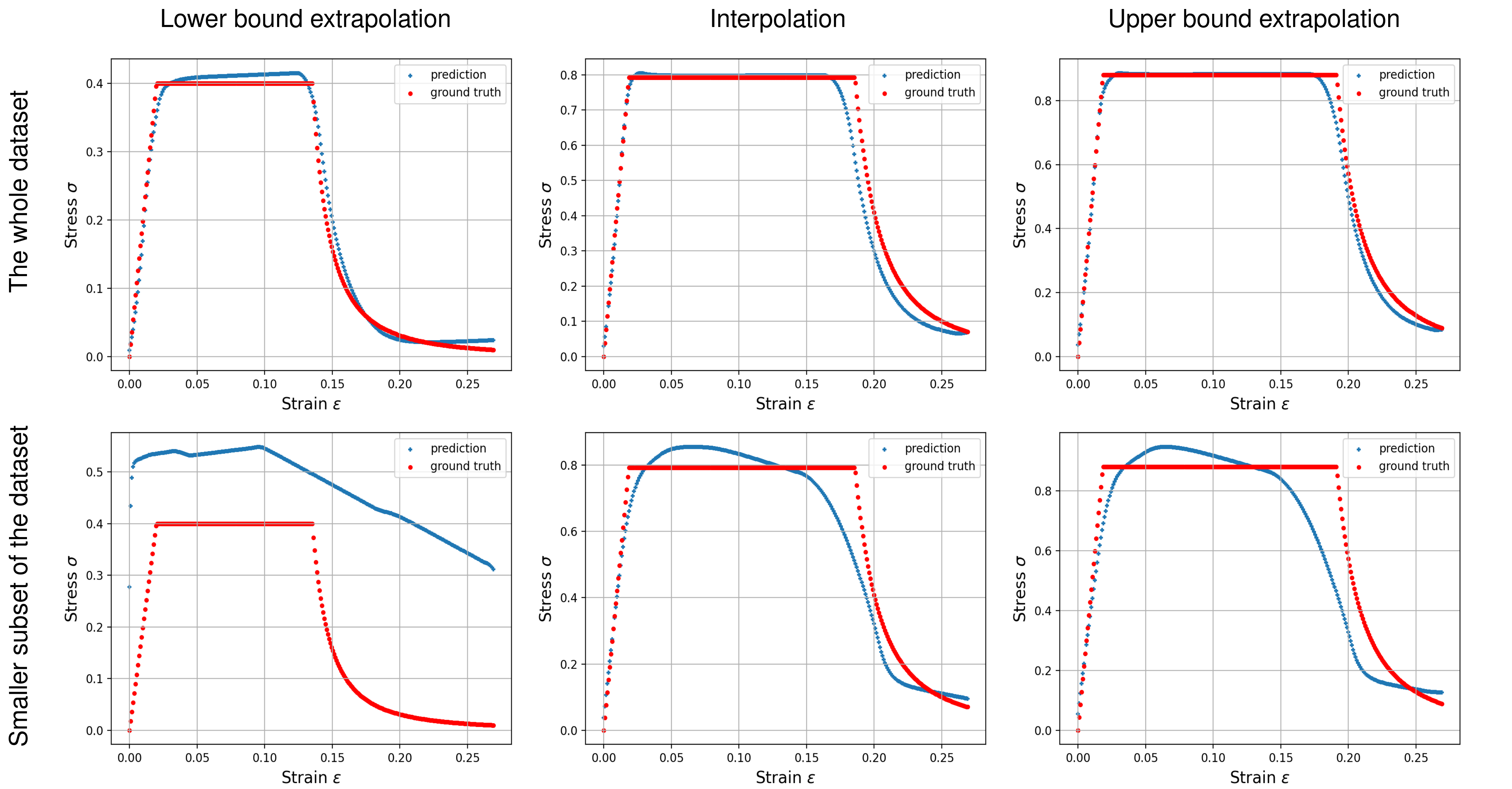} 
\caption{
Ductile fracture: Evaluation of purely data-driven, naive FFNN models employing the full (top) versus reduced datasets (bottom) for the 3 test load paths (left, middle, right column).}
\label{ductile-classical-ML}%
\end{figure}%

\begin{figure}[b]%
\centering
 \includegraphics*[width=.75\textwidth]{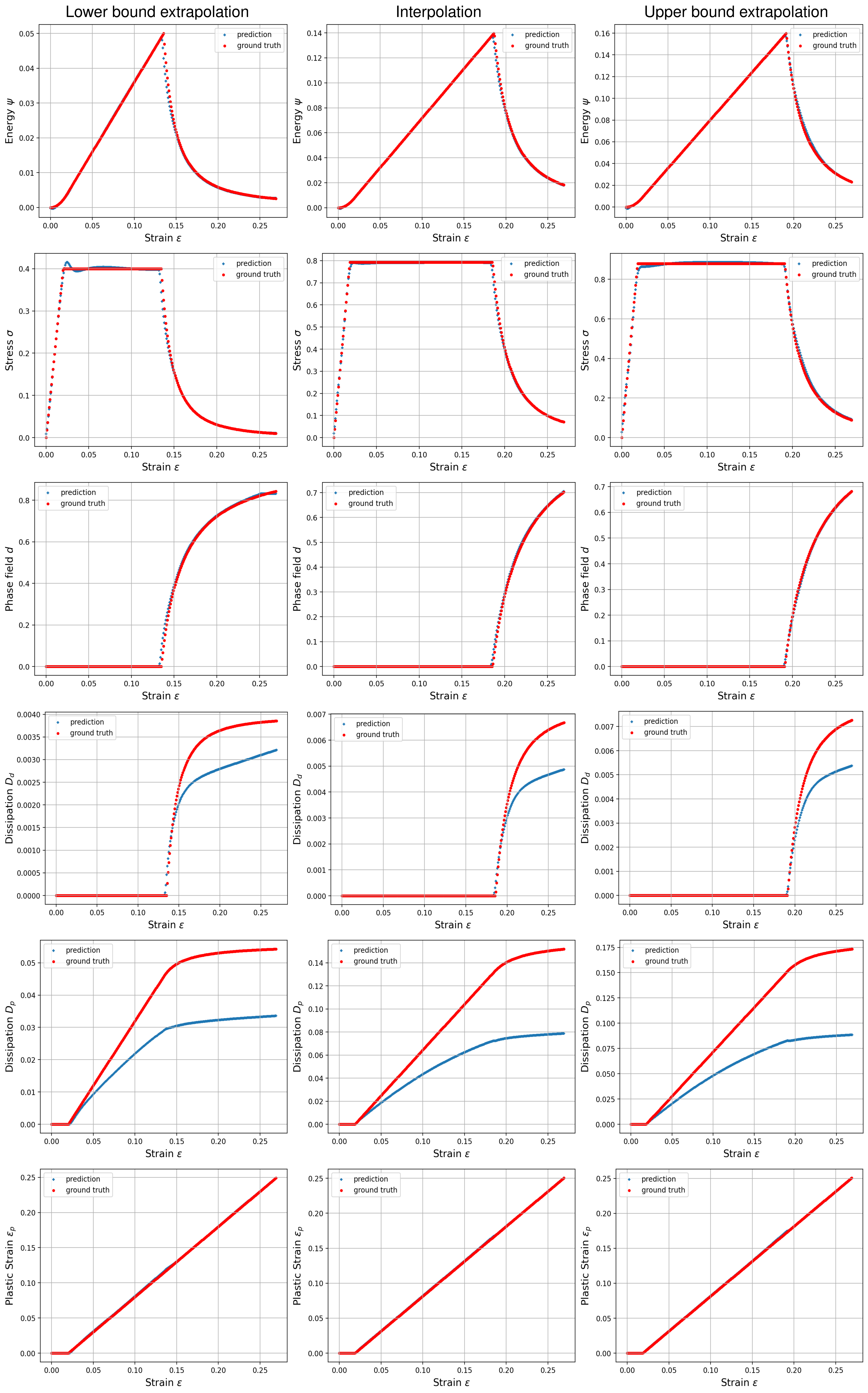
} 
\caption{
Ductile fracture: Evaluation of $\phi$ML model using the full dataset for the 3 test load paths (left, middle, right column).}
\label{ductile-pML-full}%
\end{figure}%

\begin{figure}[t]%
\centering
 \includegraphics*[width=.75\textwidth]{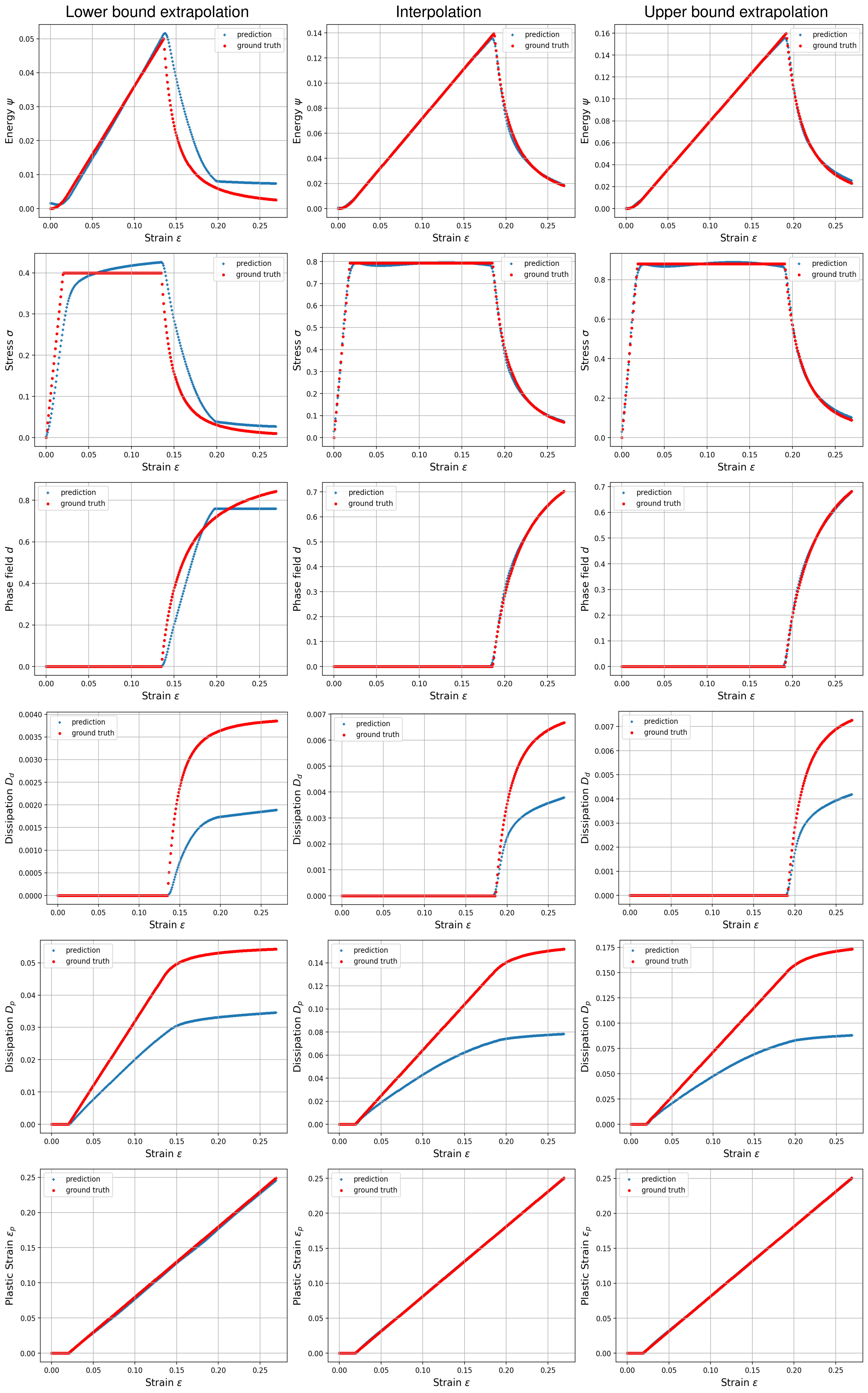
} 
\caption{
Ductile fracture: Evaluation of $\phi$ML model using the reduced dataset for the 3 test load paths (left, middle, right column).}
\label{ductile-pML-reduced}%
\end{figure}%

Table \ref{tab:comparison-brit} highlights the superior performance of the $\phi$ML model over the naive FFNN. The evaluation metrics, i.e. coefficient of determination ($R^2$) and mean absolute percentage error (MAPE), clearly illustrate this distinction. The $R^2$ value indicates how well predictions align with the ground truth, with values closer to 1 indicating higher accuracy, while MAPE quantifies the percentage error, where lower values indicate better performance. The naive FFNN struggles in extrapolation, showing very low or even negative $R^2$ values and high MAPE values in the reduced dataset, indicating poor generalization. Even in interpolation, its accuracy declines when data is limited.  In contrast, the $\phi$ML model consistently achieves $R^2$ values above 0.99 and low MAPE across both interpolation and extrapolation tasks, even with limited data, confirming its robustness, reliability, and superior generalization capabilities.

\begin{table}[!t]
    \centering
    \renewcommand{\arraystretch}{1.2}
    \begin{tabular}{l l c c c c}
        \toprule
        \textbf{Dataset} & \textbf{Region} & $R^2$ (naive FFNN) & $R^2$ ($\phi$ML) & MAPE (naive FFNN) & MAPE ($\phi$ML) \\
        \midrule
        \textbf{Full} & Lower bound extrapolation & 0.894 & 0.996 & 14.93\% & 3.18\% \\
        & Interpolation & 0.998 & 0.999 & 0.63\% & 0.41\% \\
        & Upper bound extrapolation & 0.998 & 0.999 & 3.13\% & 2.67\% \\
        \midrule
        \textbf{Reduced} & Lower bound extrapolation & -2.48 & 0.901 & 17.56\% & 14.22\% \\
        & Interpolation & 0.179 & 0.997 & 58.53\% & 7.72\% \\
        & Upper bound extrapolation & -0.182 & 0.991 & 81.09\% & 8.96\% \\
        \bottomrule
    \end{tabular}
    \caption{Comparison of $R^2$ and MAPE for naive FFNN and $\phi$ML at the brittle fracture onset position.}
    \label{tab:comparison-brit}
\end{table}

\subsection{Ductile fracture}
For the ductile failure mechanism, represented as E-P-F (Elastic-Plastic-Fracture), the material initially behaves elastically (E) up to a specific yield threshold $y_0$. Beyond this point, it undergoes plastic deformation (P). Once the deformation reaches a critical value $\ve_c = \sqrt{2\psi_c/E}$, a fracture zone develops, leading to final failure and separation of the material. The numerical results used to generate the dataset for the ML models are plotted in Figure \ref{data-generationbrit-ductile}. Hereby, a dataset of 6900 loading steps (23 load paths $\times$ 300 loading steps) is considered, to capture stress distribution, crack initiation, and propagation at each step.

Figure \ref{ductile-classical-ML} compares the performance of the data-driven, naive FFNN model for the full dataset (top row) and the reduced dataset (bottom row) in predicting stress-strain relationships for ductile fracture under three {\it test load-paths} scenarios: lower bound extrapolation, interpolation, and upper bound extrapolation. For the full dataset, the model predicts well during interpolation and captures trends in extrapolation scenarios. However, noticeable deviations occur when predicting stresses for both lower and upper extrapolation cases compared to the ground truth (red dots). With the reduced subset, the model struggles to maintain accuracy in predicting stress-strain relationships, particularly for extrapolation scenarios. These results highlight the strong dependency of data-driven models on large datasets and the challenges in generalizing beyond the training range for ductile fracture.

Figures \ref{ductile-pML-full} and \ref{ductile-pML-reduced} present the performance of the physics-based ML model in predicting ductile fracture. The model utilizes full dataset (Figure \ref{ductile-pML-full}) and reduced subset (Figure \ref{ductile-pML-reduced}), capturing key quantities, i.e. the energy $\psi$, stress $\sigma$, phase-field variable $d$, dissipation $D$, and plastic strain $\varepsilon^p$. For the full dataset, the model achieves near-perfect agreement with the ground truth (red dots) across all test load path scenarios, particularly during interpolation, and accurately captures trends in extrapolation with minimal deviations. Even with the reduced subset, the model maintains strong performance, accurately predicting interpolation and preserving the general trends in extrapolation scenarios. Minor deviations are observed in stress and dissipation predictions during extrapolation. However, the overall accuracy remains significantly higher compared to the data-driven approach (Figure \ref{ductile-classical-ML}).

To quantify this, Table \ref{tab:comparison-duc} highlights the superior performance of the $\phi$ML model in predicting ductile fracture behavior.  The naive FFNN struggles in extrapolation, with low $R^2$ values and high MAPE, especially when trained on reduced datasets. In contrast, the $\phi$ML model achieves high accuracy with $R^2$ values close to 1 and low MAPE across all cases, effectively capturing ductile fracture behavior even with reduced data.

\begin{table}[h]
    \centering
    \renewcommand{\arraystretch}{1.2}
    \begin{tabular}{l l c c c c}
        \toprule
        \textbf{Dataset} & \textbf{Region} & $R^2$ (naive FFNN) & $R^2$ ($\phi$ML) & MAPE (naive FFNN) & MAPE ($\phi$ML) \\
        \midrule
        \textbf{Full} & Lower bound extrapolation & 0.990 & 0.998 & 3.85\% & 3.73\% \\
        & Interpolation & 0.965 & 0.999 & 1.53\% & 0.02\% \\
        & Upper bound extrapolation & 0.988 & 0.998 & 0.69\% & 0.6\% \\
        \midrule
        \textbf{Reduced} & Lower bound extrapolation & -1.52 & 0.908 & 37.16\% & 6.42\% \\
        & Interpolation & 0.924 & 0.999 & 7.87\% & 0.49\% \\
        & Upper bound extrapolation & 0.833 & 0.998 & 7.6\% & 0.92\% \\
        \bottomrule
    \end{tabular}
    \caption{Comparison of $R^2$ and MAPE for naive FFNN and $\phi$ML at the ductile fracture onset position.}
    \label{tab:comparison-duc}
\end{table}

These results highlight the robustness and generalization capabilities of the physics-based ML model for ductile fracture, even with limited training data, showcasing its significant advantages over purely data-driven ML methods.

\section{Conclusion}
\label{sec6}

In this contribution, a physics-based machine learning $\phi$ML framework was developed for modeling both brittle and ductile fractures. Our framework embeds physical principles directly into the neural network architecture, unlike physics-informed neural networks that rely purely on weak enforcement of physical constraints through the loss function. This integration ensures adaptability, consistency, and a unified approach to material modeling and ML in computational fracture mechanics. The proposed approach was trained on synthetic datasets generated from finite element-based phase-field simulations, successfully capturing the homogeneous responses of brittle and ductile fractures. Detailed analyses revealed its ability to accurately predict key fracture behaviors, such as energy degradation, phase-field evolution, stress-strain relationships, and dissipation due to plasticity and fracture.

The developed $\phi$ML model addressed the limitations of classical ML models, which depend on large datasets and lack physical consistency. By leveraging its physics-integrated design, the $\phi$ML model provided accurate and reliable predictions, even with limited training data. Its strong generalization capabilities and physical consistency make it a robust tool for modeling fracture mechanics, offering a significant step toward bridging machine learning and material modeling.

While this is a first study, several open topics could be explored in future work to further enhance the proposed framework. These include the extension to complex heterogeneous anisotropic fracture behaviors, and the integration of experimental data with synthetic datasets to improve robustness and ensure greater applicability to real-world fracture problems.

\bibliographystyle{model1-num-names}
\bibliography{manuscript}

\end{document}